\documentclass[preprint,12pt]{elsarticle}

\usepackage{amsmath}
\usepackage{cases}
\usepackage{amssymb}
 \usepackage{amsthm}
 \usepackage{amsfonts}
\usepackage{multirow}
\usepackage{eqnarray}
\usepackage{diagbox} 
\usepackage{eqnarray}
\usepackage{xcolor}
\usepackage[shortlabels]{enumitem}

\usepackage{caption}
\usepackage{subcaption}

\usepackage{fancyhdr,graphicx,amsmath,amssymb}
\usepackage[ruled,vlined]{algorithm2e}
\usepackage{comment}
\usepackage{tikz}
\usepackage{pgfplots}
\pgfplotsset{compat=newest}
\allowdisplaybreaks

\newcommand{\Pstart}{{s_\text{p}}}
\newcommand{\Pend}{{e_\text{p}}}
\newcommand{\Cstart}{{s_\text{c}}}
\newcommand{\Cend}{{e_\text{c}}}
\newcommand{\CM}{{\tiny \text{\tiny CM}}}
\newcommand{\PM}{{\text{\tiny PM}}}
\newcommand{\loss}{{\text{\tiny loss}}}
\def\CMCOST#1{b_{#1}^\CM}

\def\FCMCOST#1#2{{F}_{#1}^{\CM\,#2}}

\def\CMcompCOST#1{g_{#1}^\CM}
\def\PMcost#1{c_{#1}^\PM}
\def\mobcostPMd#1{d_{#1}}
\def\mobcostPMf#1{f_{#1}}

\journal{ }
\makeatletter
\def\ps@pprintTitle{%
 \let\@oddhead\@empty
 \let\@evenhead\@empty
 \def\@oddfoot{\centerline{\thepage}}%
 \let\@evenfoot\@oddfoot}
\makeatother
\begin{document}

\begin{frontmatter}



\title{Mathematical optimization models for long-term maintenance
    scheduling of wind power systems}

\author{Quanjiang Yu and Ann-Brith Str{\"o}mberg}

\address{Department of Mathematical Sciences, Chalmers University of Technology and  
University of Gothenburg, SE-42196 Gothenburg, Sweden,\\ 
e-mail: \rm{yuqu@chalmers.se; anstr@chalmers.se}}

\begin{abstract}
\noindent
During the life of a wind farm, various types of costs arise.
A large share of the operational cost for a wind farm is due to maintenance of
    the wind turbine equipment; these costs are especially pronounced for offshore wind
    farms and they provide business opportunities in the wind energy industry. 
An effective scheduling of the maintenance activities may reduce the costs
    related to maintenance.

We combine mathematical modelling of preventive maintenance scheduling
    with corrective maintenance strategies. 
We further consider different types of contracts between the wind farm owner
    and a maintenance or insurance company, and during different phases
    of the turbines' lives and the contract periods.
Our combined preventive and corrective maintenance models are then applied
    to relevant combinations of the phases of the turbines' lives and
    the contract types.

Our case studies show that even with the same initial criteria, 
    the optimal maintenance schedules differ between different phases of time
    as well as between contract types.
One case study reveals a $40~\%$ cost reduction and a
    significantly higher production availability---$1.8~\%$ points---achieved by
    our optimization model as compared to a pure corrective maintenance strategy.
Another study shows that the number of planned preventive maintenance occasions
    for a wind farm decreases with an increasing level of an insurance contract
    regarding reimbursement of costs for broken components.


\end{abstract}
\begin{keyword}
Preventive maintenance \sep Corrective maintenance \sep Integer linear optimization \sep
    Maintenance contract \sep Wind turbine maintenance \sep Interval costs
\end{keyword}

\end{frontmatter}

\section{Introduction} \label{sec:intro}
\noindent
Global warming is an important issue nowadays from many points of view,
    since it will lead to serious global consequences for both humans and nature.
Global warming has been attributed to increased greenhouse 
    gas emission concentrations in the atmosphere through the burning of fossil fuels.
To reduce these 
    concentrations, wind power can be used to replace
    energy from fossil fuels. 
A wind farm can generate energy all day and night, as long as there is wind. 
Wind power is an abundant renewable energy provider, which produces almost no
    greenhouse gases during operation. 
According to \cite{windbenefit2020}, wind power is one of the technologies that
    produce zero greenhouse gas emissions.
Currently, the U.S.\ wind industry employs more than $100000$ workers, and
    according to \cite{windvision2020}, the wind has the potential to support more than
    $600000$ jobs in manufacturing, installation, maintenance, and supporting
    services by the year $2050$ in the U.S.
  
During the life of a wind farm, various kinds of costs arise. 
In the beginning, there are costs for, e.g., renting the premises, building permits,
    and the construction of the wind turbines. 
During the years of operation, there are maintenance costs, costs for grid
    transport of electricity, personnel costs, etc. 
When the wind turbine is outworn, costs arise for tearing it down. 
A large share of the operational cost for a wind farm is due to maintenance of
    the wind turbine equipment; these costs are especially pronounced for offshore wind
    farms and they provide business opportunities in the wind energy industry. 
An effective way to reduce the costs related to maintenance is to employ an
    improved methodology for scheduling the maintenance activities.
In this article, we investigate the planning of two types of maintenance activities,
    i.e., corrective maintenance and preventive maintenance.
    
Corrective maintenance (CM) means to maintain a component after a failure has occurred. 
Then there will be a production loss and the failure may also affect the remaining lives
    of other components. 
The cost of CM is, therefore, a combination of component costs, damage costs (for other components),
    and costs due to lost production (during the time from the failure until replacement).

Preventive maintenance (PM) means to maintain a component before it fails;
    it is a pre-scheduled maintenance activity. 
The component cost is then typically lower than the corresponding CM cost and there
    are no damage costs of other components.
Further, the downtime of the equipment is typically shorter than under CM,
    due to PM being pre-scheduled; hence, the loss of revenue from the
    production is smaller.
    
The wind turbine maintenance industry involves four types of {\sl stakeholders}:
    {\sl manufacturers}, {\sl wind farm owners}, {\sl maintenance companies},
    and {\sl insurance companies}.
The maintenance is performed by either the manufacturer, the farm owner,
    a contracted maintenance company, or a temporarily hired maintenance company;
    the interrelations of the companies involved are regulated by contracts.
We consider the following four contract types, all of which are common within
    wind energy production and maintenance.
\begin{enumerate}[label={[{\text{C}-\Roman*}]}]
\parskip -0.5mm

\item \label{full-contract}

{\sl Full service production-based maintenance contract}
    (between the wind farm owner and the manufacturer).
Maintenance is performed by the manufacturer, who covers the costs for
    replacement related to all component failures.
The contract also guarantees a minimum 'level of measured average availability'
    of the wind farm, defined as
    'available production divided by the sum of available and unavailable production'.
The manufacturer plans for PM; if a component suddenly fails the manufacturer
    performs CM of the broken component.
A variant of this type of contract includes a {\sl time-based warranty},
    which guarantees a minimum 'level of technical availability of the wind farm',
    defined as 'the total share of each given time interval during which 
    the wind farm is available for operation'.
Under this contract---since the manufacturer pays for all maintenance 
    costs---we optimize the PM scheduling on behalf of the manufacturer.

\item \label{basic-contract}
{\sl Basic insurance contract}
 (between the party who pays for the maintenance and the insurance company).
Maintenance is performed either by a maintenance company or by the wind farm owner's
    maintenance team.
When a component fails due to a 'sudden failure' there is a need
    for CM of that component.
The insurance company will reimburse the farm owner for the cost of
    the component with a discount, while the work/labour costs
    are paid by the wind farm owner.  
In practice, the stakeholders negotiate in order to classify a failure
    as a sudden failure or not; 
    in our modelling this negotiation is represented by a probability that
    'the failure is classified as being sudden'.
Since a sudden failure entails unwanted costs for both stakeholders,
    they both benefit from performing PM, where the extent depends
    on the value of the probability.
Under this contract, the wind farm owner and the insurance company share
    the components' costs during a CM activity. 
Hence, the PM scheduling is optimized on behalf of the wind farm owner.

\item \label{no-contract}
{\sl No insurance contract}.
Maintenance is performed either by a maintenance company or by the wind farm owner's
    maintenance team.
The wind farm owner covers all its costs for maintenance. 
    The wind farm owner plans for PM; if a component suddenly fails the owner asks
    a maintenance company or its own maintenance team to perform CM.
Under this contract---since the wind farm owner pays for all
    maintenance costs---the maintenance scheduling is optimized
    on behalf of the wind farm owner.
    
\item \label{maintenance-company-contract}
{\sl Maintenance contract with a maintenance company}
    (between the wind farm owner and the maintenance company).
Four main types of agreements exist:
    \begin{itemize}
    \parskip -0.5mm
    \item {\sl Call-off agreement}. 
        The simplest variant; the wind farm owner contacts the service provider in
            the event of an error. 
        The wind farm owner pays for the working hours required for
            the maintenance and for the components as the costs occur.
    \item {\sl Basic agreement}. 
        The planned service is included in an annual fee; CM
            is paid by the wind farm owner when it occurs. 
        This agreement is available both with and without remote monitoring.
    \item {\sl Full service "light"}. 
        The planned service, corrective maintenance, and spare parts,
            monitoring, as well as reporting, are included in an annual fee.
        The main components (blade/rotor, main bearing, gearbox, generator,
        nacelle, tower, and foundation) are excluded.
        Inverters and SCADA systems are either included or excluded. 
        This agreement can be with or without an availability guarantee (time- or 
        production/energy-based). 
    \item {\sl Full service}. 
        The planned service, CM (including spare parts and main components),
            monitoring, and reporting are included in the agreement. 
         Blade wear--and--tear are included in certain agreements. 
         Foundations are not included.
        This agreement comes with an availability guarantee, which is time- or
        production-based.
    \end{itemize}
We assume that related to the main component, this contract can be
    characterized as either of the contract types \ref{full-contract},
    \ref{basic-contract}, or \ref{no-contract}. 
For the call-off agreement,  the basic agreement and full service "light" agreement the optimal scheduling strategy will be
    similar to that of \ref{no-contract}.
For the full service agreement, the optimal scheduling strategy will be similar to that of
    \ref{full-contract}.
\end{enumerate}

\noindent
According to discussions with a group of wind farm owners within the Swed\-ish
    Wind Power Technology Centre (SWPTC) \cite{SWPTC},
    two different set-ups are commonly experienced. 
One set-up is to have a \ref{full-contract} contract with the manufacturer
    from the beginning till the end of life of the wind farm. 
The other set-up is to have a \ref{full-contract} contract with the manufacturer
    during the initial years of a wind farm's operating period, after which
    the wind farm owner either extends this contract (type \ref{full-contract}),
    or acquires a contract with a maintenance company
    (type \ref{maintenance-company-contract}),
    or acquires a basic insurance contract (type \ref{basic-contract}),
    or has no insurance contract (type \ref{no-contract});
    during the life of the wind farm, the owner may switch between
    the four contract types.

Regardless of the type of contract that applies during each given time period, 
    the wind turbines need to be continuously maintained throughout their lives.
From a maintenance planning perspective, the interrelations of the {\sl end of
    the planning period} considered, the {\sl end of the current contract period},
    and the {\sl end of the life} of the wind farm
    are used to define the following three {\sl phases} of time.
\begin{enumerate}[label={[{\text{P}-\roman*}]}]
\parskip -0.5mm

\item \label{normal-phase}
    The {\sl normal phase}: After the end of the planning period there is still
        a long period of time left for the contract period
        as well as for the life of the wind farm.
    The maintenance company thus cares about component failures also after
        the end of the planning period.
    We address this responsibility by imposing penalties that grow with
        the ages of the respective components that are present in the turbine
        at the end of the planning period.
    
\item \label{end-of-full-contract-phase}
    The {\sl phase close to the end of a full service production-based
        maintenance contract} (presuming a contract of type \ref{full-contract}):
    The end of the planning period coincides with the end of the contract period.
    The maintenance company needs to ensure that the wind farm is functioning at
        the contracted level of availability until the end of the planning period;
        the maintenance company is not responsible for any component failures
        after the end of the planning period.
    
\item \label{end-of-life-phase}
    The {\sl phase close to the end of life of the wind farm}
        (presuming a contract of type \ref{basic-contract} or \ref{no-contract}):
    The end of the planning period coincides with the end of life of
        the wind farm.
    Hence, it may be non-beneficial to maintain any components that fail
        (very) close to the end of the planning period.
\end{enumerate}

\noindent
We will study four principal cases that are relevant for the planning of PM
    and which are characterized by the above definitions of contracts and phases.

First, we consider the case when the wind farm and the maintenance company
    has a contract of type \ref{full-contract} and
    the maintenance company is responsible for the maintenance. 
We propose in Section \ref{sec:full-service-contract} a flexible
    mathematical optimisation model for finding optimal maintenance plans
    in the phase \ref{end-of-full-contract-phase}.
The model finds a PM schedule over the planning period such that
    the costs of performing maintenance are minimized.
The availability constraints encountered in phase 
    \ref{end-of-full-contract-phase} are modelled as the expected average
    (over the planning period) number of available turbines. 

Section \ref{sec:normal-phase} considers a wind farm owner with
    a contract of either type \ref{full-contract} or type \ref{no-contract}. 
In phase \ref{normal-phase} our model finds a PM schedule over the 
    planning period such that the costs of performing maintenance are minimized. 
The maintenance costs are supplemented by penalty costs based on the ages of 
    the component individuals in the turbine(s) at the end of the planning period.

In Section \ref{sec:end-of-life-phase} we consider a wind farm owner
    without any insurance contract, i.e., with a contract of type
    \ref{no-contract}, and in the phase \ref{end-of-life-phase}.
The model finds a PM schedule over the planning period such that the costs
    of performing maintenance are minimized;
    it also finds an optimal time, after which CM should not be performed.

Section \ref{sec:delta-insurance} considers a wind farm owner having
    a basic insurance contract, i.e., a contract of type \ref{basic-contract}.
Our model then seeks to minimize the maintenance costs, while accounting
    for a certain probability that the insurance company will pay for
    new components at CM.
The relation between the price of an insurance contract and the cost of an optimal
    maintenance plan---for the case of having no insurance contract---is
    also modelled.

The four principal cases are modelled using a framework of mathematical modelling,
    which is described in the introductory sections of the paper.
In Section \ref{sec:literature} we review relevant literature on
    (industrial) maintenance planning and contracting between
    the wind farm owner/operator and the maintenance provider.
Section \ref{sec:farm-model} describes our basic mathematical model for
    finding feasible PM schedules for a farm of wind turbines,
    the modelling of expected costs for pure CM strategies,
    costs for PM scheduling, and combinations of these;
    this results in so-called {\sl interval costs} for component replacements; 
    the modelling of costs is based on Weibull distributed failure times
    of the components.
In Section~\ref{sec:PMscheduling}, the four principal cases are described, 
    followed by the reporting in Section~\ref{sec:Case-studies} of results from 
    a set of case studies.
Finally, in Section \ref{sec:conclusion} we draw conclusions from the
    results of the case studies {and propose further research topics}.

\section{Literature review} \label{sec:literature}

\noindent
Before the 20th century, maintenance cost was considered to be a kind of inevitable cost,
    and most maintenance actions performed considered corrective maintenance. 
In 1929, the idea of preventive maintenance by daily, weekly, and general inspections was
    introduced in \cite{rose1929motorcoach}. 
Since then, maintenance methodology has grown quite fast. 
There is a broad body of literature devoted to various strategies of maintenance scheduling. 
In this article, we look at maintenance scheduling methods from different perspectives.

In the paper by \citet{yeh2006periodical}, a mathematical model to derive an optimal
    \textit{periodical PM policy} for a leased facility is developed. 
Within a lease period, any failures of the facility are rectified by minimal repairs and
    a penalty may occur to the lessor when the time required to perform a minimal repair
    exceeds a reasonable time limit. 
Further on, \citet{lee2016new} considered periodic PM policies for a deteriorating repairable system,
    and the effect of each PM action is classified into one of the three categories
    'failure rate reduction', 'decrease of deterioration speed', and 'age reduction'. 
    
While periodical PM policies consider equidistant PM occasions, 
    \citet{gustavsson2014preventive} and \citet{moghaddam2011sensitivity} study
    \textit{PM scheduling over a long time period}. 
The model PMSPIC in \cite{gustavsson2014preventive} was devised to schedule PM of
    the components of a system over a finite and discretized time horizon,
    given a common set-up cost and component costs dependent on the lengths of
    the maintenance intervals;
the model can be used for PM scheduling as well as be dynamically used
    in a setting allowing for rescheduling.
In \cite{moghaddam2011sensitivity} optimization models are developed that determine
    optimal PM schedules for repairable and maintainable systems.
It is demonstrated that higher set-up costs make simultaneous PM activities beneficial. 
The suggested models are, however, nonlinear, meaning that they are computationally hard to solve.
   
 \citet{jafari2018joint} propose a joint optimization of the maintenance policy
    and the inspection interval for a multi-unit series system. 
The authors develop a model and algorithm which are used to determine a maintenance
    policy that minimizes the maintenance cost for a multi-component system;
    in the system studied, one unit is subject to condition monitoring, while for
    the other unit, only age information is available and which has a general distribution. 
\citet{tian2014condition} develop a method that uses the condition monitoring data
    to effectively predict the remaining life of a component in a multi-component system. 
    
The papers \cite{de2001multicriteria, lisnianski2008maintenance, park2016cost, qiu2017optimal}
    are devoted to optimization issues related to \textit{different maintenance contracts}.
\citet{park2016cost} deal with the optimal maintenance policy under different warranty
    policies, considering both the warranty period and the post-warranty period. 
For the warranty period, the authors suggest a warranty cost model using a
    repair--replacement warranty policy considering repair as well as
    failure times.  
\citet{qiu2017optimal} consider optimization of the maintenance costs under
    performance-based contracts.
The paper investigates an optimal maintenance policy for inspected systems that are
    subject to both soft and hard failures. 
According to \citet{de2001multicriteria}, the main parameters of maintenance contracts
    are downtime and maintenance costs. 
\citet{lisnianski2008maintenance} consider an aging system, in which the maintenance
    is performed by an external maintenance team. 
Different kinds of contracts between the two parties are considered,
    with special attention paid to downtime costs
    (i.e., the loss of revenue due to production stops). 
The authors suggest a model based on a piece-wise constant approximation of 
    the increasing failure rate function.


\section{Mathematical modelling of maintenance planning for a farm of wind turbines} \label{sec:farm-model}
\noindent
We consider the maintenance planning for a wind farm comprising $m$ wind turbines
    (indexed by $i \in \mathcal{I} := \{ 1, \ldots, m \}$)
    each of which has $n$ (identical) component types
     (indexed by $j \in \mathcal{J} := \{ 1, \ldots, n \}$).     

During the usage of a wind turbine its components will degrade and---if not
    repaired or replaced---one or more components will eventually fail;
    then, CM will have to be performed in the form of component replacement.
We will model the failure rates such that the survival function of a component of
    type $j \in \{ 1, \ldots, n \}$ comes from a Weibull distribution with
    scale parameter $\alpha_j >0$ and shape parameter $\beta_j >0$. 

\begin{figure}[htb]
    \centering
    \begin{tikzpicture}[scale=0.7]
\foreach \x in {0,1,2,5,6,7,10,11,14,15,16}
    \draw [shift={(\x,0)}, color=black] (0pt,5pt) -- (0pt,-5pt);
    
\node[above] at (0,0.2) {\footnotesize $0$};
\node[above] at (1,0.2) {\footnotesize $1$};
\node[above] at (2,0.2) {\footnotesize $2$};
\node[above] at (3.5,0.2) {\footnotesize $\cdots$};
\node[above] at (6,0.2) {\footnotesize $\Pstart$};
\node[above] at (7,0.2) {\footnotesize $\Pstart\!\!+\!\!1$};
\node[above] at (8.5,0.2) {\footnotesize $\cdots$};
\node[above] at (10,0.2) {\footnotesize $\Pend\!\!-\!\!1$};
\node[above] at (11,0.2) {\footnotesize $\Pend$};
\node[above] at (12.5,0.2) {\footnotesize $\cdots$};
\node[above] at (15,0.2) {\footnotesize $T$};

\node[above] at (0.5,-1.2) {\footnotesize $1$};
\node[above] at (1.5,-1.2) {\footnotesize $2$};
\node[above] at (3.5,-1.2) {\footnotesize $\cdots$};
\node[above] at (5.5,-1.2) {\footnotesize $\Pstart$};
\node[above] at (6.5,-1.2) {\footnotesize $\Pstart\!\!+\!\!1$};
\node[above] at (8.5,-1.2) {\footnotesize $\cdots$};
\node[above] at (10.5,-1.2) {\footnotesize $\Pend$};
\node[above] at (12.5,-1.2) {\footnotesize $\cdots$};
\node[above] at (14.5,-1.2) {\footnotesize $T$};

\node[below] at (0.5,0) {\footnotesize $\underbrace{\hphantom{\;}}$};
\node[below] at (1.5,0) {\footnotesize $\underbrace{\hphantom{\;}}$};
\node[below] at (5.5,0) {\footnotesize $\underbrace{\hphantom{\;}}$};
\node[below] at (6.5,0) {\footnotesize $\underbrace{\hphantom{\;}}$};
\node[below] at (10.5,0) {\footnotesize $\underbrace{\hphantom{\;}}$};
\node[below] at (14.5,0) {\footnotesize $\underbrace{\hphantom{\;}}$};

\draw[-] (-0.6,0) -- (16.6,0);
\end{tikzpicture}
    \caption{Illustration of the notation for the continuous and normalized time
        (above the time line) and the corresponding time steps (below the time line)
        used in the models.}
    \label{fig:intervals-and-steps}
\end{figure}
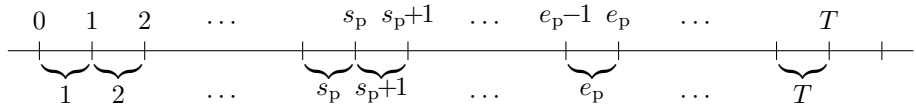

Our models are defined over a time period---denoted by the {\sl continuous time interval}
    $[0,T]$---such that a (new) wind turbine starts production at time $0$,
    while its life ends at time $T$.
Part of our modelling employs a finite set of {\sl discrete time steps}%
\footnote{$T$ is assumed to be an integer and the time scale is assumed to be normalized,
            such that each time step has the length of one (1) time unit.},
    indexed by $\{ 1, \ldots, T \}$, where time step $t \in \{ 1, \ldots, T \}$ represents
    the time interval $[t-1, t] \subset [0, T]$.
A maintenance contract (abbreviated as $\text{c}$) period covers a time interval
    $[\Cstart,\Cend] \subseteq [0,T]$, such that $\Cend \geq \Cstart+1$, i.e.,
    a time period covering the time steps indexed by 
    $\{ \Cstart+1, \ldots, \Cend \} \subseteq \{ 1, \ldots, T \}$.
Analogously, a maintenance planning (abbreviated as $\text{p}$) period covers a time interval
    $[\Pstart,\Pend] \subseteq [0,T]$, such that  $\Pend \geq \Pstart+1$, i.e.,
    a time period covering the time steps indexed by
    $\{ \Pstart+1, \ldots, \Pend \} \subseteq \{ 1, \ldots, T \}$.
Figure~\ref{fig:intervals-and-steps} illustrates the time steps and the
    corresponding continuous time intervals. 

\subsection{A basic mathematical model defining a feasible PM schedule}
\label{sec:PM-planning}

The mathematical model \eqref{eq:PMSPIC-farm}, below, is an adaption to our settings
    of the model of the preventive maintenance scheduling problem with
    so-called {\sl interval costs},%
\footnote{The concept of 'interval cost' is defined as the cost of maintaining a component
            as a function of the time interval from the component's previous maintenance occasion.}
    as developed in \cite{gustavsson2014preventive}.
It considers the planning period $(\Pstart,\Pend]$, indexed by $\{ \Pstart+1, \ldots, \Pend \}$,
    during which PM actions can be scheduled.

The variables of the model are defined as follows.
For each wind turbine $i \in \mathcal{I}$ and component type $j \in \mathcal{J}$,
    and any two indices $u \in \{ \Pstart, \ldots, \Pend \}$ and
    $t \in \{ u+1, \ldots, \Pend+1 \}$, $x^{ij}_{ut} = 1$ if PM is planned at time steps
    $u$ and $t$, but no PM is planned in-between these;%
\footnote{For any component $(i,j)$ we adopt the following interpretation conventions: 
        For any $t \in \{ \Pstart+1, \ldots, \Pend \}$
            $x^{ij}_{\Pstart t} = 1$ ($x^{ij}_{t,\Pend+1} = 1$) means that the first (last) PM is planned at time step $t$,
        while $x^{ij}_{\Pstart,\Pend+1} = 1$ means that no PM at all is planned.}
 otherwise $x^{ij}_{ut} = 0$.
Whenever $x^{ij}_{ut} = 1$, we say that a PM interval for component $(i,j)$ starts
    at time step $u$ and ends at time step $t$.
For each wind turbine $i \in \mathcal{I}$ and any two indices $u \in \{ \Pstart, \ldots, \Pend-1 \}$ and
    $t \in \{ \Pstart+1, \ldots, \Pend \}$, the variable $y^{i}_{ut}=1$ if PM is planned at time steps
    $u$ and $t$, but no PM is planned in-between these 
    for at least one of its component types at time $t$; otherwise $y^{i}_{ut} = 0$;
    see \eqref{tu}.
Further, for any two indices $u \in \{ \Pstart, \ldots, \Pend-1 \}$ and $t \in \{ \Pstart+1, \ldots, \Pend \}$, the variable
    $z_{ut} = 1$ if PM is planned at time steps
    $u$ and $t$, but no PM is planned in-between these for any of the components in the farm at time $t$;
    otherwise $z_{ut} = 0$; see \eqref{fa}.
\begin{subequations} \label{eq:PMSPIC-farm}
\begin{align}
 && \sum_{t=\Pstart+1}^{\Pend+1} x^{ij}_{\Pstart t} & = 1,\ && 
        i \in \mathcal{I},\ j \in \mathcal{J}, \label{one} \\
    && \sum_{u=\Pstart}^{t-1} x^{ij}_{ut} & = \sum_{v=t+1}^{\Pend+1} x^{ij}_{tv}, && 
        i \in \mathcal{I},\ j \in \mathcal{J},\ && t=\Pstart+1, \ldots, \Pend,\  
        \label{ba} \\
    && y_{ut}^i & \geq x^{ij}_{ut}, && 
        i \in \mathcal{I},\ j \in \mathcal{J},\ && u = \Pstart, \ldots, t-1,
        \label{tu} \\ 
    &&  &   &&  && t = \Pstart+1, \ldots \Pend,\ 
        \nonumber \\
    && z_{ut} & \ge y_{ut}^i, && 
       i \in \mathcal{I},\ && u = \Pstart, \ldots, t-1,
       \label{fa} \\ 
    &&  &   &&  && t = \Pstart+1, \ldots \Pend, \nonumber \\
    && x^ {ij}_{ut} & \in \{0,1\},\ &&
       i \in \mathcal{I},\ j \in \mathcal{J},\ && u = \Pstart, \ldots, t-1, \\ 
    &&  &   &&  && t = \Pstart+1, \ldots \Pend+1,\ \nonumber \\
    && y^{i}_{ut}, z_{ut} & \in \{0,1\},\ &&
       i \in \mathcal{I},\ && u = \Pstart, \ldots, t-1, \\ 
    &&  &   &&  && t = \Pstart+1, \ldots \Pend.\ \nonumber
\end{align} 
\end{subequations}
Due to the constraints \eqref{one}, for each component $(i,j)$ either the first PM
    is scheduled in one of the time steps $\{ \Pstart+1, \ldots, \Pend \}$,
    or no PM is scheduled (i.e., $x^{ij}_{\Pstart,\Pend+1} = 1$). 
The constraints \eqref{ba} make sure that, for each component $(i,j)$, 
    the end of any PM interval is the start of the next PM interval.%
\footnote{The constraints \eqref{one}--\eqref{ba} are equivalent to so-called
            flow balance constraints; see \cite{fulkerson1966flow}.}

Consider a point which is feasible subject to the constraints \eqref{eq:PMSPIC-farm},
    and which is such that for all $(i,j)$ and all $u \in \{ \Pstart, \ldots, t-1 \}$
    and $t \in \{ \Pstart+1, \ldots, \Pend+1 \}$, $x^{ij}_{ut} = 0$, except that
    $x^{ij}_{\Pstart,\Pend+1} = 1$, and such that for all $i$ and all
    $t \in \{ \Pstart+1, \ldots, \Pend \}$, $y_{ut}^i = z_{ut} = 0$.
This point corresponds to planning no PM activities during the planning period
    $[\Pstart, \Pend]$ and is referred to as a {\sl pure CM plan};
    the corresponding situation will be modelled further in Section~\ref{noPM}.

\subsection{Modelling of costs for CM strategies and PM scheduling} \label{noPM}

We next define the parameters used for the types of costs incurred by
    different maintenance activities, and describe the calculation of
    the expected maintenance costs of a CM strategy, and how these interact
    with a planned PM schedule.
We also introduce a number of parameters to be used as building blocks of the models
    to be defined in the following sections.

\subsubsection{Cost and time parameters} \label{sec:cost-time-param}

For each component type $j \in \mathcal{J}$, we define the following time-independent
    parameters representing cost for maintenance activities:
    $\CMcompCOST{j}$ denotes the cost of a new component of type $j$ at a CM action; 
    $\CMCOST{j}$ denotes the cost of
    a new component (i.e., $\CMcompCOST{j}$) plus the logistics costs
    (i.e., transport, set-up, and work costs for the crane, and manpower costs
    for the replacement) of a CM action;
$\PMcost{j}$ denotes the cost
    of a new component minus the value (i.e., the expected sales revenue)
    of the replaced component, plus the cost (e.g., work cost for the crane and
        manpower cost for the replacement) of the PM action.

The costs and revenues that are common to all component types,
    but dependent on the time step $t \in \{ 1, \ldots, T \}$, are defined as follows:
$\mobcostPMd{t}$ denotes the cost attributed to a specific wind turbine,
    incurred by any maintenance activity at time $t$;
$\mobcostPMf{t}$ denotes the costs shared by the wind turbines of a farm (e.g., 
    crane transport and set-up costs) incurred by any maintenance activity at time $t$;
$r_t  \geq 0$ denotes the revenue (in terms of the monetary value of the energy produced)
    created by a fully functioning turbine during time step $t \in \{ 1, \ldots, T \}$,
    i.e., during the continuous time period $[t-1, t] \subset [0, T]$. 

The amount of time required for performing a CM (of component $j \in \mathcal{J}$) and
    a PM action are denoted by $\gamma_j^\CM > 0$ and $\gamma^\PM > 0$, respectively.

\subsubsection{Revenue function and revenue loss} \label{sec:revenue-function}
The total revenue created by a functioning wind turbine during the time period
    $[v+\lambda,T]$, where $v \in \{ 0, \ldots, T-1 \}$ and $ \lambda \in [0, 1]$,
    is expressed by the {\sl revenue function}
    $R: [0,T+ \bar{\gamma}] \mapsto \mathbb{R}_+$, where
    $\bar{\gamma} := \max \{ \gamma^\PM; \max_{j \in \mathcal{J}} \{ \gamma_j^\CM \} \}$,
    and which is defined by
\begin{subequations} \label{eq:R(T)}
\begin{align}
    R(u) & := 0, && u \in [T, T + \bar{\gamma} ]; \label{eq:R(T):T} \\
    R(v + \lambda) & := (1-\lambda) r_{v+1} + R(v+1), && \lambda \in [0, 1], \,
        v \in \{ T-1, \ldots, 0 \}. \label{eq:R(T):0-T-1}
\end{align}
\end{subequations}
From the definition \eqref{eq:R(T)} follows that the revenue function $R$ is
    piece-wise affine, continuous, and non-increasing on the time interval 
    $[0,T + \bar{\gamma}]$.
    
The {\sl revenue loss} caused by a production stop during the time interval between
    the two (continuous) time points $u \in [0,T]$ and $u+w \in [u,T + \bar{\gamma}]$
    is then given by
\begin{align} \label{eq:Rloss}
    R^\loss(u, w) & := R(u) - R(u+w).
\end{align}

Utilizing the definitions in Section~\ref{sec:cost-time-param} we conclude
    that the revenue loss caused by a CM action of component $j$ at the (continuous)
    {\sl time} $t \in [0, T]$ equals $R^\loss(t, \gamma_j^\CM)$,
    while the revenue loss caused by a PM action of any component at
    {\sl time step} $t \in \{ 1, \ldots, T \}$ equals $R^\loss(t, \gamma^\PM)$.

\subsubsection{Failure times and induced CM costs}

For each wind turbine $i \in \mathcal{I}$ and component type $j \in \mathcal{J}$
    a number of component individuals will consecutively be put to work in the turbine.
Let $\tau_u^{ij} \in [ 0, u ]$ denote the age of the current component individual at
    (the end of) time step $u \in \{ 1, \ldots, T \}$.
Further, for $k=1,2, \ldots$, let $L_k^{ij} >0$  denote the life of the $k$'th consecutive
    component individual in use in the turbine, where $k=1$ corresponds to the current
    (i.e., at time step $u$) individual.
The $k$'th failure time for a component of type $j$ in turbine $i$, after the (current)
    time step $u$ (i.e., after the continuous time $u$), is then expressed as
\begin{align*}
    U_{uk}^{ij} & := u -\tau_u^{ij} + \sum_{\ell=1}^k L_\ell^{ij}, && k = 1, 2, \ldots.
\end{align*}
We also define $U_{u0}^{ij} := u - \tau_u^{ij}$ and
    assume independent, identically Weibull distributed individual component lives
    $L_k^{ij} \sim \text{W}(\alpha_j,\beta_j)$, $k=1, 2, \ldots$ .
Specifically, the life distribution of the current individual (i.e., $k=1$) is conditioned
    on its age $\tau_u^{ij}$ at time step $u$.
Hence, it holds that
\begin{align*}
    P(L_1^{ij} \geq v \,|\, L_1^{ij} \geq \tau^{ij}_u) & = \left\{
    \begin{array}{ll}
        1, & v \in [0, \tau^{ij}_u], \\
        \exp \bigg[ \Big(\frac{\tau^{ij}_u}{\alpha_j} \Big)^{\beta_j} 
            - \Big(\frac{v}{\alpha_j} \Big)^{\beta_j} \bigg], \qquad \quad
            & v \geq \tau^{ij}_u,
    \end{array} \right.
\end{align*}
and
\begin{align*}
    P(L_k^{ij} \geq v) & = \exp \bigg[ {- \Big(\frac{v}{\alpha_j} \Big)^{\beta_j}} \bigg],
        \qquad v \geq 0, \quad k=2,3, \ldots .
\end{align*}

Any failure of a component type $j \in \mathcal{J}$ in a turbine
    $i \in \mathcal{I}$ at time $t \in [0,T]$ induces a cost of a CM action
    which amounts to
\begin{align} \label{eq:CM-action-cost}
    \CMCOST{j} + R^\loss(t, \gamma_j^\CM).
\end{align}

\subsubsection{Expected cost of a pure CM strategy} \label{sec:pureCMcost}

We denote by $\FCMCOST{ut}{ij}$ the expected cost associated with the maintenance
    of a component type $j \in \mathcal{J}$ in a turbine $i \in \mathcal{I}$
    during a time interval $[u,t]$ under a \textit{pure CM strategy},
    which is such that no PM is scheduled during the time interval considered.
It is defined as the expectation of the sum of component costs and revenue losses, as
\begin{align} \label{eq:FMCOSTutij}
    \FCMCOST{ut}{ij} := \mathbb{E} \left[ \sum_{k=1}^\infty 
        1_{\{ U_{u k}^{ij} \leq t \}} \Big( \CMCOST{j} 
        + R^\loss(U_{u k}^{ij}, \gamma_j^\CM) \Big) \right].
\end{align}
Algorithm~\ref{alg:sumE} describes how to compute $\FCMCOST{ut}{ij}$.
The expected sum of component costs and revenue losses associated with CM actions
    of a wind farm under a pure CM strategy during the time interval $[\Pstart,\Pend]$
    is then defined as
    $\FCMCOST{\Pstart\Pend}{} := \sum_{i \in \mathcal{I}} \sum_{j \in \mathcal{J}}         
        \FCMCOST{\Pstart\Pend}{ij}$.

\begin{algorithm}[htb]
\SetAlgoLined
{\small
\KwResult{$\FCMCOST{ut}{ij}$}
\KwIn{$i \in \mathcal{I}$, $j \in \mathcal{J}$, $u \in \{ 1, \ldots, T-1 \}$,
    $t \in \{ u+1, \ldots, T \}$, $\tau_u^{ij} \in [0, u]$, $M \gg 1$,
    $\CMCOST{j}$, $R:\ [0,T] \mapsto \mathbb{R}_+$} 
 \For{$\ell=1, \ldots, M$}{
    $k := 0$ ; $U_{u0}^{ij} := u - \tau_u^{ij}$ ; $a_\ell := 0$ \;
    \While{$U_{uk}^{ij} \leq t \,$ }{
        $k := k + 1$ \;
        Generate $L_{k}^{ij}$ based on the survival function \;
        $U_{uk}^{ij} := U_{u,k-1}^{ij} + L_k^{ij}$ \;
        $a_\ell := a_\ell + \CMCOST{j} + R^\loss(U_{uk}^{ij}, \gamma_j^\CM)$ \;
    }
}
 $\FCMCOST{ut}{ij} := M^{-1} \sum_{\ell=1}^{M} a_\ell$ \;
}
 \caption{\small Numerical computation of expected CM costs $\FCMCOST{ut}{ij}$,
    defined in \eqref{eq:FMCOSTutij}, for a component type $j$ in a turbine
    $i$ during the time interval $[u,t]$}
    \label{alg:sumE}
\end{algorithm}

\subsubsection{Failure times and PM cost savings} \label{sec:PM-cost-savings}

In the case when PM is scheduled but components fail before their plan\-ned PM,
    a need for CM actions arises.

Consider a schedule for the PM of a wind farm over a specific time interval,
    and let $u$ and $t > u$ denote two time steps within that interval
    at each of which there is a planned PM action, but no PM is planned
    in-between these two time steps, for a component type
    $j \in \mathcal{J}$ in a turbine $i \in \mathcal{I}$.

If a component type $j \in \mathcal{J}$ fails at a time $\sigma \in [u,t]$,
    then a CM action must be performed at time $\sigma$
    (inducing a CM cost according to \eqref{eq:CM-action-cost}).
A rescheduling of the planned PM due to this failure then leads to a possible reduction of
    the costs $\PMcost{j}$ associated with the planned PM action for the component type $j$
    and turbine $i$ at time step $t$.
We denote this reduction by the \textit{PM cost savings function} 
    $\phi_{ut}^\PM: [u,t] \mapsto [0,1]$.
For $\sigma=u$, meaning that the component individual would fail immediately after
    it is put to use in the turbine and that a CM action is indeed performed at
    time step $u$, the state of this component type is retained and the current
    PM schedule remains optimal; hence, it holds that $\phi_{ut}^\PM(u) = 0$. 
For $\sigma=t$, the rescheduling would amount to removing the scheduled PM action
    for component type $j$ at time step $t$, and replacing it by a CM action;
    then, all the associated PM costs are saved and it holds that $\phi_{ut}^\PM(t) = 1$.

It is reasonable to assume that the PM cost savings function is non-decreasing with
    the failure time $\sigma \in [u,t]$,
    and we will further assume that the function $\phi_{ut}^\PM$ is linear.
This results in the expression
\begin{align*}
    \phi_{ut}^\PM(\sigma) := 
        \frac{\sigma-u}{t-u}, \qquad \sigma \in [u,t].
\end{align*}
Considering the case of $n \geq 1$ failures 
    at times $\sigma_k$, $k=1, \ldots, n$, such that
    $u \leq \sigma_1 \leq \sigma_2 \leq \ldots \leq \sigma_n \leq t$,
    the PM cost savings function is generalised as
\begin{align} \label{eq:multiple-saving-new}
    \phi_{ut}^\PM(\sigma_1, \ldots, \sigma_n) &:= 
        \frac{\sigma_1-u}{t-u} + \sum_{k=2}^n \frac{\sigma_k-\sigma_{k-1}}{t-u} 
    = \frac{\sigma_n-u}{t-u}.
\end{align}
From \eqref{eq:multiple-saving-new} follows that only the time $\sigma_n$ of the
    last failure within the time interval $[u,t]$, for component type $j \in \mathcal{J}$
    in a turbine $i \in \mathcal{I}$, will affect the PM cost savings function.

Recalling that $U_{uk}^{ij}$ denotes the time of the $k$'th failure for a component
    type $j$ in a turbine $i$, after the time step $u$, we conclude that the time
    of the last failure of the component can be expressed as
\begin{subequations}
\begin{align}
    \max \left\{\, u \, ; \, \max_{k \geq 0} \left\{\, U_{uk}^{ij} \,\big|\,
        U_{uk}^{ij} \leq t \,\right\} \,\right\}.
\end{align}
We then regard the effect of component failures on the savings of the costs
    $\mobcostPMd{t}+ R^\loss(t,\gamma^\PM)$ and $\mobcostPMf{t}$, on the wind turbine and wind farm levels, respectively. 
Reasoning analogously as above, the time $\sigma_n$ of the last failure in the
    time interval $[u,t]$ of any component in a turbine $i \in \mathcal{I}$
    is expressed as 
\begin{align}
    \max \left\{ u ;
    \max_{j \in \mathcal{J}} \left\{ \max_{k \geq 0} \left\{ U_{uk}^{ij} \,\big|\, U_{uk}^{ij} \leq t \right\} \right\} \right\},
\end{align}
    while the time $\sigma_n$ of the last failure in the time interval $[u,t]$
    of any component in any turbine in the wind farm, is expressed as
\begin{align}
    \max \left\{ u ;
    \max_{i \in \mathcal{I}, j \in \mathcal{J}} \left\{ \max_{k \geq 0} 
        \left\{ U_{uk}^{ij} \,\big|\, U_{uk}^{ij} \leq t \right\} \right\} \right\}.
\end{align}
\end{subequations}

\subsubsection{Failure-free time intervals}

For each wind turbine $i \in \mathcal{I}$, each component type $j \in \mathcal{J}$, 
    each time step $u \in \{ 1, \ldots, T-1 \}$ of the previous PM of the component,
    and each $t \in \{ u+1, \ldots, T \}$, we can now express the
    \textit{expected proportion of the failure-free, right-most subset} of
    the time interval $[u,t]$ as
\begin{subequations} \label{eq:p-ab}
\begin{align} \label{eq:pijut-ab}
    p_{ut}^{ij} := 
        \frac{1}{t-u} \cdot \left( t - \mathbb{E} \left[ \max \left\{ u ; \max_{k \geq 0}
        \left\{ U_{uk}^{ij} \,\big|\, U_{uk}^{ij} \leq t \right\} \right\} \right] \right).
\end{align}
These proportions can be numerically computed according to Algorithm~\ref{alg:Emax}.

For each wind turbine $i \in \mathcal{I}$,
    each time step $u \in \{ 1, \ldots, T-1 \}$ of the previous PM occasion for the turbine,
    and each $t \in \{ u+1, \ldots, T \}$, the corresponding failure-free proportion
    of the time interval $[u,t]$ equals 
\begin{align} \label{eq:piut-ab}
    p_{ut}^i := \frac{1}{t-u} \cdot \left( t - \mathbb{E} \left[ \max \left\{ u ;
    \max_{j \in \mathcal{J}} \left\{ \max_{k \geq 0} \left\{ U_{uk}^{ij} \,\big|\,
    U_{uk}^{ij} \leq t \right\} \right\} \right\} \right] \right).
\end{align}

On the wind farm level, for each time step $u \in \{ 1, \ldots, T-1 \}$ of
    the previous PM occasion for the farm, and each $t \in \{ u+1, \ldots, T \}$,
    the corresponding failure-free proportion of the time interval $[u,t]$ equals 
\begin{align} \label{eq:put-ab}
    p_{ut} := \frac{1}{t-u} \cdot \left( t - \mathbb{E} \left[
    \max \left\{ u ;
    \max_{i \in \mathcal{I}, j \in \mathcal{J}} \left\{ \max_{k \geq 0} 
    \left\{ U_{uk}^{ij} \,\big|\, U_{uk}^{ij} \leq t \right\} \right\} \right\}
    \right] \right).
\end{align}
\end{subequations}

\begin{algorithm}[htb]
\SetAlgoLined
{\small
\KwResult{$p_{ut}^{ij}$}
\KwIn{$i \!\in\! \mathcal{I}$, $j \!\in\! \mathcal{J}$, $u \!\in\! \{ 1, \ldots, T\!-\!1 \}$,
    $t \!\in\! \{ u\!+\!1, \ldots, T \}$, $\tau_u^{ij} \!\in\! [0, u]$, $M \!\gg\! 1$}
\For{$\ell=1, \ldots, M$}{
    $k := 0$ ;
    $U_{u0}^{ij} := u -\tau_u^{ij}$ \;
    \Repeat{$U_{uk}^{ij} > t$}{
        $k := k+1$ \; 
        Generate $L_{k}^{ij}$ based on the survival function \;
        $U_{uk}^{ij} := U_{u,k-1}^{ij} + L_k^{ij}$ \; 
        $a_\ell := U_{u,k-1}^{ij}$}
    \If{$k=1$}{$a_\ell := u$}
    }
$p_{ut}^{ij} := (t-u)^{-1} \big( t - M^{-1} \sum_{\ell=1}^{M} a_\ell \big)$
}
 \caption{\small Numerical computation of $p_{ut}^{ij}$, defined in \eqref{eq:pijut-ab},
    for a component type $j$ in a turbine $i$ during the time interval $[u,t]$}
    \label{alg:Emax}
\end{algorithm}

\subsubsection{Interval costs}

We define the {\sl interval cost} of a PM action at time step $t$ on a component type $j$
    of turbine $i$, when the previous maintenance of the same component type was performed
    at time step $u \leq t-1$,
    as the sum of the cost of the scheduled PM action at time $t$,
    the revenue loss caused by the PM action at time $t$, and the expected CM costs
    for the component type $j$ of turbine $i$ over the time interval $[u,t]$.
For the case when a CM action {\sl is not} accompanied by a rescheduling of the PM
    schedule (see Section \ref{sec:pureCMcost}), the interval cost over
    the time interval $[u,t]$ is thus defined as
\begin{align*} 
    \FCMCOST{ut}{ij} + \PMcost{j},
\end{align*}   
    where $\FCMCOST{ut}{ij}$ is defined in \eqref{eq:FMCOSTutij}.
For the case when a CM action {\sl is} accompanied by a rescheduling of the PM schedule,
    part of the cost $\PMcost{j}$ can be saved,
    and the interval cost over the time interval $[u,t]$ is defined as
\begin{align*} 
    \FCMCOST{ut}{ij} + \PMcost{j} p_{ut}^{ij},
\end{align*}   
    where $p_{ut}^{ij}$ is defined in \eqref{eq:pijut-ab}.

\section{Modelling of PM scheduling under different phases and contracts}
\label{sec:PMscheduling}
\noindent
The four combinations of contracts and phases that are relevant for the planning of PM 
    of wind farms are described in the following subsections.

\subsection{Full service production-based contract} \label{sec:full-service-contract}

Here, we consider the phase close to the end of a full service production-based
    maintenance contract, i.e., phase \ref{end-of-full-contract-phase} and
    contract type \ref{full-contract}.
The planning period is then defined as the time interval $[\Pstart,\Cend]$.

Hence, for $t = \Pstart+1, \ldots, \Cend+1$, $u = \Pstart, \ldots, t-1$, $i \in \mathcal{I}$,
    and $j \in \mathcal{J}$, we define the sum of expected component costs and revenue losses as%
\footnote{Since there is no PM after the planning period, we define $p_{u, \Cend+1}^{ij} := 0$.}
\begin{subequations} \label{eq:C-ut-ijny}
    \begin{numcases}{C_{ut}^{ij} :=} 
        \FCMCOST{ut}{ij} + \PMcost{j} p_{ut}^{ij}, & $t \in \{ \Pstart\!+\!1, \ldots, \Cend \}$, 
        \label{eq:C-ut-ijny:a} \\
        \FCMCOST{u,\Cend}{ij}, & $t = \Cend\!+\!1$, \label{eq:C-ut-ijny:b}
    \end{numcases}
\end{subequations}
    where $\FCMCOST{ut}{ij}$ and $p_{ut}^{ij}$, are defined in \eqref{eq:FMCOSTutij}
    and \eqref{eq:pijut-ab}, respectively. 
Letting $p_{ut}^i$ and $p_{ut}$ be defined by \eqref{eq:piut-ab} and \eqref{eq:put-ab},
    respectively, for any pair $(s,e)$ of indices such that $s \geq 0$ and 
    $e \in \{ s+1, \ldots, T \}$, we define the function $G^\PM_{se}$ of PM costs
    over the time interval $[s,e]$, on the turbine and farm levels as
\begin{align} \label{eq:PMcostTurbFarm}
    G^\PM_{se}(\mathbf{y},\mathbf{z}) := \sum_{u=s}^{e-1}     
        \sum_{t=u+1}^{e} \left( \Big( \mobcostPMd{t}+ R^\loss(t,\gamma^\PM) \Big) \sum_{i \in \mathcal{I}} 
        p_{ut}^i y^{i}_{ut} + \mobcostPMf{t} p_{ut} z_{ut} \right).
\end{align}
Letting $C^{ij}_{ut}$ be defined by  \eqref{eq:C-ut-ijny}, an optimal PM plan
    is then obtained by minimizing the sum of expected costs for CM and the costs for
    the PM occasions during the time period $[\Pstart, \Cend]$, according to
\begin{subequations} \label{model:full-contract}
\begin{align} \label{eq:F-ab}
    \mathop{\text{minimise }}_{\mathbf{x},\mathbf{y},\mathbf{z}} & \sum_{u=\Pstart}^{\Cend} \sum_{t=u+1}^{\Cend+1} 
   \sum_{i \in \mathcal{I}} \sum_{j \in \mathcal{J}} C^{ij}_{ut} x^{ij}_{ut} 
    + G^\PM_{\Pstart \Cend}(\mathbf{y},\mathbf{z}) \qquad \\
    \text{ subject to} & \text{ the constraints in } \eqref{eq:PMSPIC-farm}.
\end{align}
\end{subequations}

\subsubsection{A guaranteed level of availability}

A {\sl full service production-based maintenance contract} \ref{full-contract}
    guarantees that the {\sl measured average availability} of the wind farm
    is not below a certain level. 
We will model this guarantee as the requirement that the 
    {\sl expected average availability} is above a certain level.

\paragraph{The expected average availability being above a certain level}
\noindent
For a planning period $[\Pstart+1, \Cend]$, the total revenue from a fully functioning
    wind turbine equals $R^\loss(\Pstart, \Cend \!-\! \Pstart)$
    (see Section \ref{sec:revenue-function}).
The requirement that the average availability is above a certain level can
    then be expressed as putting an upper limit 
    $\epsilon \cdot R^\loss(\Pstart, \Cend \!-\! \Pstart)$ on the loss of revenue
    due to maintenance actions being performed, where $\epsilon \in (0,1]$.

The upper limit on the revenue loss, due to production downtime
    during the time period $[\Pstart+1, \Cend]$
    is thus expressed by the inequality constraint 
    (cf.\ the expressions in~\eqref{eq:C-ut-ijny})
\begin{align} \label{D}
    \nonumber
    & \sum_{u=\Pstart}^{\Cend} \sum_{t=u+1}^{\Cend+1} \sum_{i \in \mathcal{I}} 
        \sum_{j \in \mathcal{J}} \!\Bigg(\! \mathbb{E} \!\left[ \sum_{k=1}^{\infty}
        1_{\{ U_{uk}^{ij} \leq \min \{ t, \Cend \} \} } R^\loss(U_{uk}^{ij},\gamma_j^\CM) 
        \right] \!\Bigg) x^{ij}_{ut} \\ 
    & \qquad \qquad \qquad + \sum_{u=\Pstart}^{\Cend-1}
    \sum_{t=u+1}^{\Cend} R^\loss(t,\gamma^\PM) \sum_{i \in \mathcal{I}} p_{ut}^{i} y_{ut}^{i}
      \leq \epsilon R^\loss(\Pstart, \Cend \!-\! \Pstart).
\end{align}
An optimal PM plan is then defined by a solution to the
    integer linear optimisation problem to
\begin{subequations} \label{model:full-contract-guaranteed-availability}
\begin{align} 
    \mathop{\text{minimise }}_{\mathbf{x},\mathbf{y},\mathbf{z}} &  
    \text{the objective in } \eqref{eq:F-ab}, \label{model:full-contract-guaranteed-availability-obj}\\ 
    \text{subject to } & \text{the constraints } \eqref{eq:PMSPIC-farm} \text{ and } \eqref{D}.
    \label{model:full-contract-guaranteed-availability-constr} 
\end{align}
\end{subequations}

\paragraph{Full service time-based maintenance contract}

\noindent
The variant of the full service production-based maintenance contract
    \ref{full-contract} including a {\sl time-based warranty} guarantees
    a minimum level of technical availability of the wind farm,
    defined as 'the total share of each given time interval during which 
    the wind turbines are available for operation'.
For the planning period $[\Pstart+1,\Cend]$, this guarantee is then modelled as
    the time loss due to downtime not exceeding the upper limit
    $\epsilon [\Cend - \Pstart]$, where $\epsilon \in (0,1]$.

The interval cost is then altered to a time-based measure, defined as
\begin{align*} 
    \gamma_j^\CM \, \mathbb{E} \Bigg[ \sum_{k=1}^{\infty}
        1_{\{ U_{uk}^{ij} \leq \min \{ t, \Cend \} \} } \Bigg]
    = \gamma_j^\CM \, \mathbb{E} \left[ \max_{k \geq 0}
        \left\{\, k \,\big|\, U_{uk}^{ij} \leq \min \{ t, \Cend \} \,\right\} \right],
\end{align*}
    which can be computed by a slightly modified version of Algorithm~\ref{alg:Emax}.

The price of a full service time-based maintenance contract is thus given by
    the minimum value of the objective in \eqref{eq:F-ab} 
    being further constrained by the upper limit on the time loss due to
    downtime, as
\begin{subequations} \label{model:full-time-guaranteed-availability}
\begin{align} 
    & \mathop{\text{minimise }}_{\mathbf{x},\mathbf{y},\mathbf{z}} \text{the objective in } \eqref{eq:F-ab}, \\
    & \text{subject to the constraints } \eqref{eq:PMSPIC-farm} \text{ and} \\
    & \nonumber \sum_{u=\Pstart}^{\Cend} \sum_{t=u+1}^{\Cend+1} \sum_{i \in \mathcal{I}}
        \sum_{j \in \mathcal{J}} \gamma_j^\CM \ \mathbb{E} \Big[ \max_{k \geq 0}
        \left\{ k \,\big|\, U_{uk}^{ij} \leq \min \{ t, \Cend \} \right\} \Big] x^{ij}_{ut} \\
    & \qquad \qquad \qquad \qquad \qquad \quad + \sum_{u=\Pstart}^{\Cend-1} \sum_{t=u+1}^{\Cend} 
        \gamma^\PM \sum_{i \in \mathcal{I}} p_{ut}^{i} y_{ut}^{i} \; \leq \;
        \epsilon ( \Cend - \Pstart ). 
\end{align}
\end{subequations}

\subsection{Optimizing the PM plan in the normal phase} \label{sec:normal-phase}

We consider optimizing the maintenance costs in the normal phase,
    i.e., \ref{normal-phase}, either with a
    {\sl full service production-based maintenance contract} or with
    {\sl no insurance contract}, i.e., either \ref{full-contract} or \ref{no-contract}.
For these two types of contracts, the optimization models are equal, but the 
    {\sl stakeholder} who pays the cost for the maintenance (both PM and CM) differs.
The time period considered is defined by the start and the end of the planning period,
    i.e., the time interval $[\Pstart,\Pend]$.

As modelled in Section~\ref{sec:full-service-contract},
    {\sl close to the end of a full service production-based maintenance contract},
    i.e., in the phase~\ref{end-of-full-contract-phase}, the maintenance company
    only needs to make sure that the components will survive until the end of the
    contract period.
In the {\sl normal phase}, i.e., the phase \ref{normal-phase}, however, at the end of
    the planning period---since the turbines are aimed to produce electricity 
    also after the planning period---the ages of the components in the turbines
    make sense for the wind farm owner.\footnote{In phase~\ref{end-of-full-contract-phase}---since the end of the planning period
    equals the end of the contract---a component failure after the planning period
    is not the maintenance company's responsibility. 
    In phase~\ref{normal-phase}---since the end of the planning period is far from
    both the end of the contract and the end of life---the party who pays
    for the maintenance cost cares about the wind turbine be functioning also
    after the end of the planning period. 
    Hence, the ages of components at the end of the planning period matters,
    and there should be a compensating fee, which increases with the ages of 
    the components.}

We model the aim of retaining a functioning wind turbine by defining penalties for 
    used components "left" in the turbine at the end of the planning period. 
For $t \in \{ \Pstart+1, \ldots, \Pend+1 \}$, $u \in \{ \Pstart, \ldots, t-1 \}$,
    $i \in \mathcal{I}$, and $j \in \mathcal{J}$, 
    the interval costs (cf.\ \eqref{eq:C-ut-ijny}) are thus altered as
\begin{subequations} \label{eq:interval-cost-altered}
    \begin{numcases}{c_{ut}^{ij} :=}
            \FCMCOST{ut}{ij} + \PMcost{j} p_{ut}^{ij},
         & $t \in \{ \Pstart+1, \ldots, \Pend - \lfloor \gamma^\PM \rfloor \}$, 
         \label{eq:interval-cost-altered:a} \\
        \FCMCOST{uT}{ij} - \FCMCOST{t-1,T}{ij}, 
        & $t \in \{ \Pend+1 - \lfloor \gamma^\PM \rfloor, \ldots, \Pend+1 \}$,
        \label{eq:interval-cost-altered:b}
    \end{numcases}
\end{subequations}
    where $\FCMCOST{ut}{ij}$ and $p_{ut}^{ij}$ are defined in
    \eqref{eq:FMCOSTutij}  and \eqref{eq:pijut-ab}, respectively,
    and $T$ denotes the end of life of the wind turbine. 
The optimal PM plan is then a solution to the problem to
\begin{subequations}
\begin{align} \label{model:normal}
    \mathop{\text{minimise }}_{\mathbf{x},\mathbf{y},\mathbf{z}} & \sum_{u=\Pstart}^{\Pend} \sum_{t=u+1}^{\Pend+1}
        \sum_{i \in \mathcal{I}} \sum_{j \in \mathcal{J}} c^{ij}_{ut} x^{ij}_{ut} 
        + G^\PM_{\Pstart\Pend}(\mathbf{y},\mathbf{z}), \\
    \text{subject to} & \text{ the constraints } \eqref{eq:PMSPIC-farm},
\end{align}
\end{subequations}
    where $G^\PM_{\Pstart \Pend}(\mathbf{y},\mathbf{z})$ is defined in \eqref{eq:PMcostTurbFarm}.

\subsection{Lifetime planning without an insurance} \label{sec:end-of-life-phase}

We next consider the time close to the end of life of the turbine(s) and 
    a wind farm owner without any insurance contract, i.e.,
    a contract of type \ref{no-contract} during the phase \ref{end-of-life-phase}.
Then, at a component failure it may not be beneficial to perform CM,
    as explained next.
As defined in \eqref{eq:R(T)}, for any $v \in [0,T]$, $R(v)$ equals the revenue
    created during the time interval $[v,T]$.
Therefore, if a component $j$ that fails at time $v \in [0,T-\gamma_j^\CM]$ is replaced,
    the maximum revenue created after the replacement equals $R(v+\gamma_j^\CM)$.
Hence, performing CM of a component type $j$ at time $v$ is profitable
    only if this revenue exceeds the component and logistics cost for CM, i.e., 
    if $R(v+\gamma_j^\CM) > \CMCOST{j}+R^\loss(U_{uk}^{ij}, \gamma_j^\CM)$ holds.
We define the cost of a failure at time $v$ as the minimum of
$\CMCOST{j}$ and $R(v+\gamma_j^\CM)$, and the corresponding modified interval
    costs by\footnote{Since there is no PM at time $T+1$, we define $p_{u,T+1}^{ij} := 0$.}
 (cf.\ the expressions in \eqref{eq:C-ut-ijny})
\begin{subequations} \label{eq:interval-cost-vC} 
    \begin{numcases}{\check{C}_{ut}^{ij} \!:=\!}
        \mathbb{E} \!\left[ {\displaystyle \sum_{k=1}^{\infty} 1_{\{ U_{uk}^{ij} \leq t \}} }
            \min \!\Big\{ \CMCOST{j} \!+\! R^\loss(U_{uk}^{ij}, \gamma_j^\CM) ; \;
            R(U_{uk}^{ij} \!+\! \gamma_j^\CM) \Big\} \right] + \PMcost{j} p_{ut}^{ij}, \nonumber \\
        \qquad \qquad \qquad \qquad \qquad \qquad \qquad
            t \in \{ \Pstart\!+\!1, \ldots, T \}, \label{eq:interval-cost-vC:a} \vspace{4mm} \\
        \mathbb{E} \!\left[ {\displaystyle \sum_{k=1}^{\infty} 1_{\{ U_{uk}^{ij} \leq T \}} }
            \min \!\Big\{ \CMCOST{j} \!+\! R^\loss(U_{uk}^{ij}, \gamma_j^\CM) ; \; R(U_{uk}^{ij} \!+\!
            \gamma_j^\CM ) \Big\} \right]\!, \nonumber \\
        \qquad \qquad \qquad \qquad \qquad \qquad \qquad \qquad \qquad t = T\!+\!1.
            \label{eq:interval-cost-vC:b}   
    \end{numcases}
\end{subequations}
The optimal PM schedule is then the solution to the optimization problem
\begin{subequations} \label{eq:hatFmodel}
\begin{align} 
    \mathop{\text{minimise }}_{\mathbf{x},\mathbf{y},\mathbf{z}} &
        \sum_{u=\Pstart}^{T} \sum_{t=u+1}^{T+1} 
        \sum_{i \in \mathcal{I}} \sum_{j \in \mathcal{J}} \check{C}^{ij}_{ut} x^{ij}_{ut} + G^\PM_{\Pstart T}(\mathbf{y},\mathbf{z}) \\
    \text{subject to } & \text{the constraints } \eqref{eq:PMSPIC-farm},
\end{align}
\end{subequations}
    where $G^\PM_{\Pstart T}(\mathbf{y},\mathbf{z})$ is defined in
    \eqref{eq:PMcostTurbFarm}.
Then, if there is a component failure at time $v \in [0, T-\gamma_j^\CM]$, 
    the component should be replaced only if it holds that
    $R(v+\gamma_j^\CM) \geq \CMCOST{j}+R^\loss(v,\gamma_j^\CM)$.

\subsection{The price of a $\delta$-insurance contract} \label{sec:delta-insurance}

Considering the case when the wind farm owner has a basic insurance contract, i.e.,
    \ref{basic-contract}, phase \ref{end-of-full-contract-phase}. 
We let $\delta \in [0,1]$ denote the probability that, at a component failure, 
    the insurance company will pay for a new component. 
We call such a contract, a {\sl $\delta$-insurance contract}.
Denoting the cost of a new component of type $j$---at a CM action---by $\CMcompCOST{j}$ 
    yields the modified interval cost (cf.\ the expressions in \eqref{eq:C-ut-ijny})
\begin{subequations} \label{eq:Cdot} 
    \begin{numcases}{\dot{C}_{ut}^{ij}(\delta) \!:=\!}
        \mathbb{E} \left[ {\displaystyle \sum_{k=1}^{\infty} 1_{\{ U_{uk}^{ij} \leq t \}}} 
        \!\!\left( \CMCOST{j} \!-\! \delta \CMcompCOST{j} \!+\! R^\loss(U_{uk}^{ij},\gamma_j^\CM)
        \right) \right] + \PMcost{j} p_{ut}^{ij}, \nonumber \\
    \qquad \qquad \qquad \qquad \qquad \qquad \qquad \qquad
        t \in \{ \Pstart\!+\!1, \ldots, \Cend \}, \label{eq:Cdot:a} \vspace{2mm} \\
    \mathbb{E} \!\left[ {\displaystyle \sum_{k=1}^{\infty} 1_{\{ U_{uk}^{ij} \leq \Cend \}}}
        \!\!\left( \CMCOST{j} \!-\! \delta \CMcompCOST{j} \!+\! R^\loss(U_{uk}^{ij},\gamma_j^\CM)
         \right) \right]\!, \quad t \!=\! \Cend\!+\!1. \qquad \quad \label{eq:Cdot:b}
    \end{numcases}
\end{subequations}
The optimal PM plan is then obtained as a solution to the following
    minimization problem, where $F^*(\delta)$ denotes the value of a
    $\delta$-insurance contract:
\begin{subequations} \label{model:delta-insurance}
\begin{align}
    F^*(\delta) := \mathop{\text{minimum }}_{\mathbf{x},\mathbf{y},\mathbf{z}} & 
    \sum_{u=\Pstart}^{\Cend} \sum_{t=u+1}^{\Cend+1} 
        \sum_{i \in \mathcal{I}} \sum_{j \in \mathcal{J}} \dot{C}^{ij}_{ut}(\delta)
        x^{ij}_{ut} + G^\PM_{\Pstart \Cend}(\mathbf{y},\mathbf{z}) \\
    \text{subject to } & \text{the constraints } \eqref{eq:PMSPIC-farm},
\end{align}
\end{subequations}
    where $G^\PM_{\Pstart \Cend}(\mathbf{y},\mathbf{z})$ is defined in
    \eqref{eq:PMcostTurbFarm}.
For the case of a $0$-insurance contract, the wind farm owner always has to pay for
    the components at CM; the value of such a contract is denoted by $F^*(0)$.
For any value of $\delta \in [0, 1]$ the difference
    \[ F^*(0) - F^*(\delta) \]
    can thus be viewed as the maximum price that a wind farm owner should pay for a $\delta$-insurance contract. 
A corresponding altering of the costs $\CMCOST{j}$ to 
    $\CMCOST{j} - \delta \CMcompCOST{j}$ in the models
    \eqref{model:full-time-guaranteed-availability} and \eqref{eq:hatFmodel}, 
    can be employed for adjusting these models with respect to different phases of time.

\section{Case studies} \label{sec:Case-studies}

\noindent 
All computational tests were performed on an Intel 2.40 GHz dual core Windows PC
    with 16 GB RAM. 
The mathematical optimization models are implemented in AMPL IDE
    (version 12.1; see \cite{ampl}),
    the parameters of the model are calculated by Matlab 
    (version R2019b; see \cite{matlab}), 
    and the optimization models are solved using CPLEX
    (version 20.1; see \cite{cplex}).

All the case studies performed consider a wind farm with ten wind turbines. 
The life of a wind turbine is assumed to be $20$ years, which is the typical case
    for onshore wind farms; see \cite{ZIEGLER20181261}. 
We studied a case in which all the components were initially 'as good as new'. 
The first time step (month) of the planning period is January.

The data for component lives (shape and scale parameters of the Weibull distribution),
    as well as the costs for PM and CM listed in Table~\ref{tab:Task1_parameters},
    are from the article by \citet[Tables 1 and~2]{tian2011condition}.
We chose the cost parameter $\mobcostPMf{t} = \$ 50000$ 
    (see \citet[Table~2]{tian2011condition}) and assume that the cost that
    can be saved on the turbine level is the downtime cost
    (i.e., the revenue loss due to downtime), i.e., $\mobcostPMd{t} := \$ 0$.

\begin{table}[!ht]
\newcommand{\tabincell}[2]{\footnotesize \begin{tabular}{@{}#1@{}}#2\end{tabular}}
\centering
{\small
\begin{tabular}{@{}l@{\qquad}c@{\qquad}c@{}c@{\qquad}c@{\qquad}c@{}}
    \hline
    \tabincell{c}{Component type, $j$} & 
    \multicolumn{2}{c@{}}{\tabincell{c}{Replacement costs}} &&
    \multicolumn{2}{c@{}}{\tabincell{c}{Weibull parameters}} \\
    \tabincell{c}{{}\\{}} & 
    \tabincell{c}{corrective\\$\CMCOST{j} \ {[\$ 1000]}$} &
    \tabincell{c}{preventive\\$\PMcost{j} \ {[\$ 1000]}$} &&
    \tabincell{c}{shape\\$\beta^j$} &
    \tabincell{c@{}}{{scale}\\{$\alpha^j$ [months]}} \\ 
    \hline
    {\footnotesize Rotor} & {\footnotesize 162} & {\footnotesize 28} && {\footnotesize 3} & {\footnotesize 100} \\ 
    {\footnotesize Main bearing} & {\footnotesize 110} & {\footnotesize 15} && {\footnotesize 2} & {\footnotesize 125} \\
    {\footnotesize Gearbox} & {\footnotesize 202} & {\footnotesize 38} && {\footnotesize 3} & {\footnotesize 80} \\
    {\footnotesize Generator} & {\footnotesize 150} & {\footnotesize 25} && {\footnotesize 2} & {\footnotesize 110} \\
    \hline
    {\footnotesize Downtime cost, $R^\loss$} & {\footnotesize 15--30} & {\footnotesize 2.5--5} && & \\
    \hline
\end{tabular}
}
  \caption{Parameters of the survival function and CM and PM costs for each
  component, and the intervals for downtime costs due to PM and CM}
  \label{tab:Task1_parameters}
\end{table}

The monthly revenue $r_t$ is calculated monthly productions multiplied
    by monthly prices and then averaged over three years. 
In February (September) the monthly revenue is the lowest (highest),
    around $\$15000$ ($\$30000$). 
These numbers are gathered from a specific wind farm,
    which had a lot of problems with icing,
    which in turn mean that the production during the winter season
    was not as high as predicted from the average wind speeds.
The downtime caused by PM and CM activities are $\gamma^\PM=\frac{1}{6}$ and
    $\gamma_j^\CM=1$ (months), $j \in \mathcal{J}$, respectively; 
    these numbers are consulted with experts within the SWPTC. 
From these figures, the revenue losses as defined in~\eqref{eq:Rloss},
    caused by the downtime due to PM and CM activities in month $t$,
    are determined by 
    $R^\loss(t,\gamma^\PM) \in [15 \gamma^\PM, 30 \gamma^\PM] 
        = [\$2500, \$5000]$ and
    $R^\loss(t,\gamma^\CM) \in [15 \gamma^\CM, 30 \gamma^\CM] 
        = [\$15000, \$30000]$, respectively
        (see Table~\ref{tab:Task1_parameters}). 

At time step (month) $t=0$, all ten turbines are assumed to be new,
    possessing identical statuses as well as conditions.
This means that in an optimal maintenance plan all ten turbines will possess
    identical maintenance schedules.
Hence, for this case, the models can be simplified to consider only one turbine,
    with all cost entities multiplied by ten,
    except the mobilization costs $\mobcostPMf{t}$ which are shared by the turbines.
This yields a computationally more tractable model.

\subsection{Parameter analysis}
\label{parameter_analysis}
Our case studies are introduced by a comparison of the interval cost
    of replacing the gearbox, for different replacement intervals $[u,t]$, 
    where $t \geq u+1$ and $u \geq \Pstart$, and during different phases of time.
The types of interval costs compared are 
    (a)~$C_{ut}^{ij}$ and $c_{ut}^{ij}$ as defined in~\eqref{eq:C-ut-ijny:a}
        and~\eqref{eq:interval-cost-altered:a}, respectively;
    (b)~$c_{ut}^{ij}$ as defined in~\eqref{eq:interval-cost-altered:b};
    (c)~$C_{ut}^{ij}$ as defined in~\eqref{eq:C-ut-ijny:b};
    (d)~$\check{C}_{ut}^{ij}$ as defined in~\eqref{eq:interval-cost-vC:b}.
For all four types, {\sl each interval starts at a PM occasion} at time step
    $u \geq \Pstart$. 
For type~(a), {\sl each interval ends at a PM occasion} at a time step $t \geq u+1$
    in either phase \ref{normal-phase}, \ref{end-of-full-contract-phase}, 
    or~\ref{end-of-life-phase}, while 
    for the types~(b)--(d), {\sl each interval ends at the end of the planning period, without a PM occasion},
    i.e., 
    for type~(b), $t=\Pend+1$ in the phase \ref{normal-phase}; 
    for type~(c), $t=\Cend+1$ in the phase \ref{end-of-full-contract-phase}
        and a contract type \ref{full-contract};
    for type~(d), $t=T+1$ in the phase~\ref{end-of-life-phase}
        and a contract type \ref{basic-contract} or~\ref{no-contract}.
    

\begin{figure}[h]
\centering
\begin{subfigure}[b]{0.47\textwidth}
    \centering
    \includegraphics[width=\textwidth]{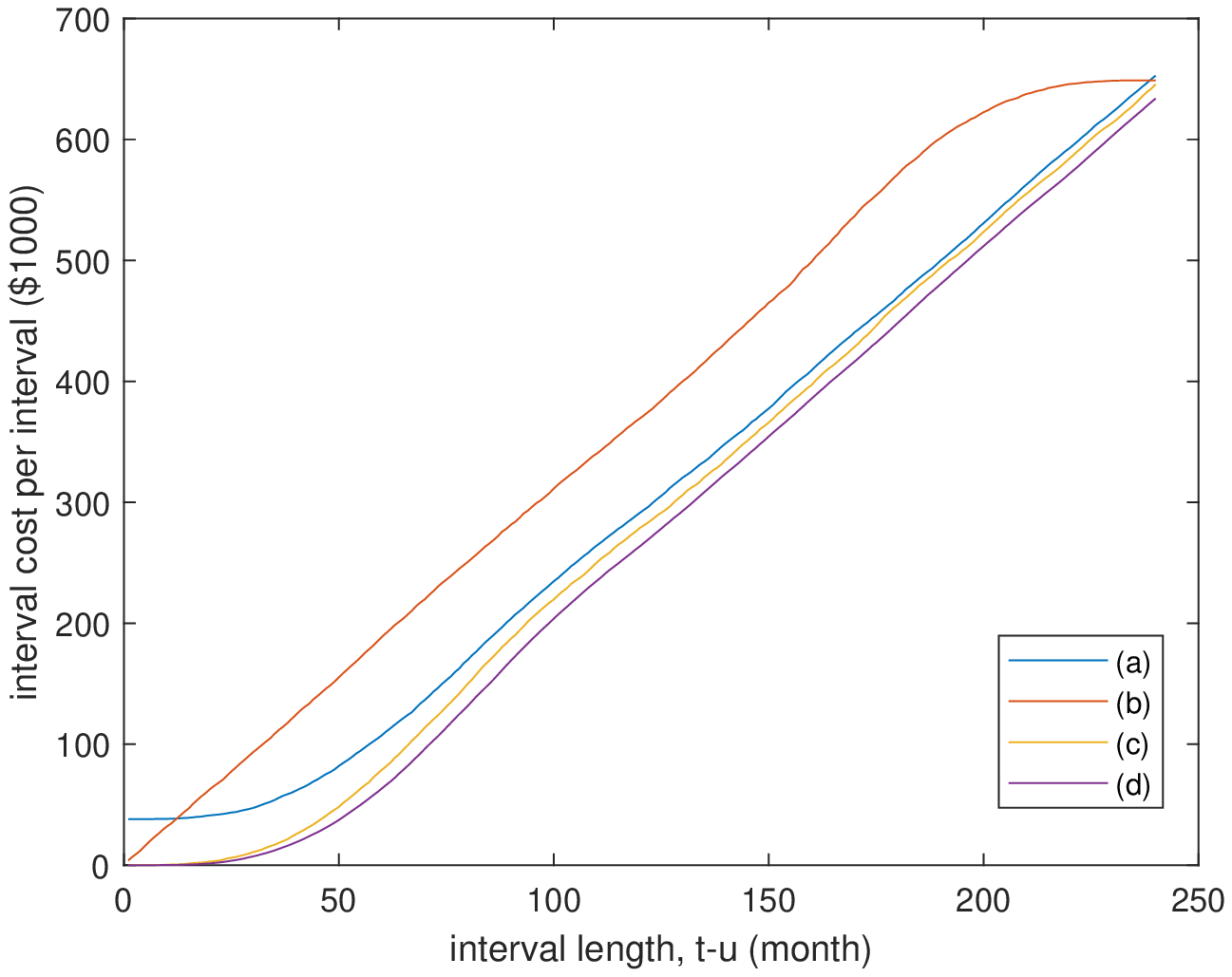}
    \caption{\small Interval cost per interval as a function of the interval length $t-u$}
    \label{fig:compare_interval_cost}
\end{subfigure}
\hfill
\begin{subfigure}[b]{0.47\textwidth}
\centering
    \includegraphics[width=\textwidth]{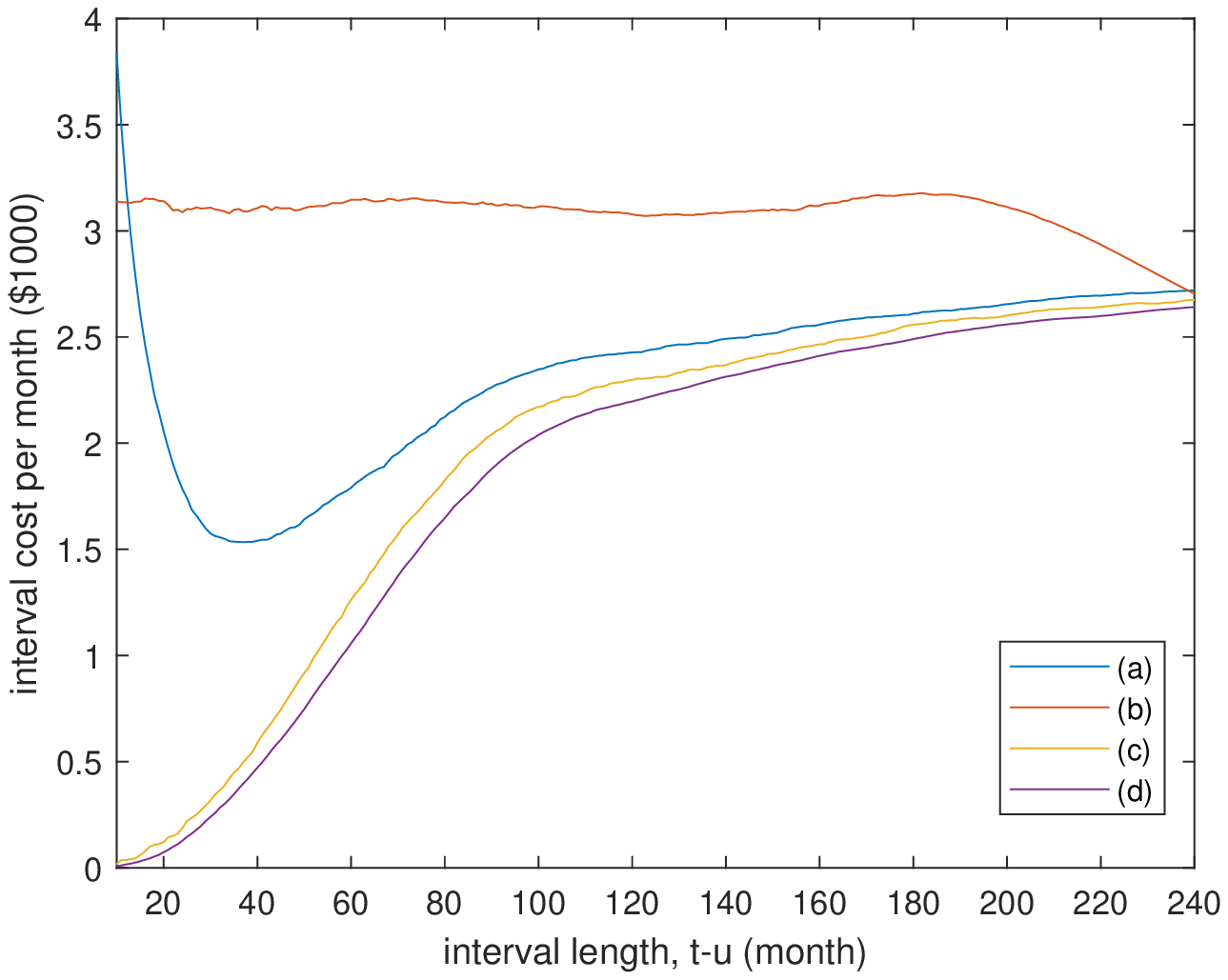}
    \caption{\small Interval cost per month as a function of the interval length $t-u$}
    \label{fig:cost_per_month}
\end{subfigure}
\caption{Comparison of interval costs for the gearbox during different phases of time: 
    Each interval starts at a PM occasion and ends at the
    (a) next PM occasion; 
    (b)--(d) end of the planning period, which is
    (b) long before the end of the contract as well as the end of life;
    (c) at the end of the contract;
    (d) at the end of life.}
\label{fig:intervalcosts}
\end{figure}
    

For interval cost type~(b), at the end of the interval, there is still a long contract period
    as well as component life left.
Figure~\ref{fig:intervalcosts}(\subref{fig:compare_interval_cost}) reveals that    
    the interval cost type~(b) is significantly higher than those of types~(c) and~(d),
    which is due to the gearbox being old with a high risk of break-down and that
    interval cost type~(b) takes into account the age of the gearbox also after
    time step $t$.
For long intervals (i.e., when $t-u$ is close to $T$), the wind farm's age
    is close to its end of life.
Hence, $c_{ut}^{ij}$ is close to $C_{ut}^{ij}$, i.e., the interval costs
    of type (b) tend to those of type (c) for long intervals. 
When $t$ is at the end of the planning period, it holds that
    $\check{C}_{ut}^{ij} < C_{ut}^{ij}$ for all $u < t$, which is due to
    the fact that if the gearbox breaks too close to the end of its life,
    maintenance will not be beneficial; as a result of not performing
    maintenance, production will be lost.
For any $u < t < \Cend$, i.e., when $t$ denotes a PM occasion, the inequalities
    $C_{ut}^{ij} > C_{\Cend-t+u,\Cend}^{ij}$ and 
    $c_{ut}^{ij} > C_{\Cend-t+u,\Cend}^{ij}$ hold, because PM costs need to be paid
    at PM occasions but not at the end of the contract period
    (i.e., at time step $\Cend$).
However, with an increasing probability of a breakdown before
    the scheduled PM activity, the difference between these costs decreases.


Figure~\ref{fig:intervalcosts}(\subref{fig:cost_per_month}) shows that for the 
    gearbox and interval cost type~(a), i.e., the interval ends at a PM occasion,
    the maintenance cost per month is the lowest at the interval length $t-u = 38$ (months).
For interval cost type~(b), when $t-u$ approaches $T$, the probability 
    of a breakdown during the time interval $[t, T]$ tends to $0$;
    hence, the monthly maintenance cost tends to the interval cost type (c). 
For interval costs type~(c) and type~(d), the monthly maintenance costs
    sometimes increase faster than the corresponding cost type~(b); 
    hence, it may be beneficial to employ a longer planning period for
    type~(b) than for type~(c) or type~(d).

\subsection{Case study 1: A ten years full service production-based maintenance contract} \label{sec:10years_contract}

We consider a full service production-based maintenance contract,
    i.e., \ref{full-contract}, over ten years, i.e., $120$ months, in 
    phase~\ref{end-of-full-contract-phase}.
Our aim is to find the most beneficial scheduling, in terms of in which 
    months the maintenance (i.e., replacement) of the different components
    should be performed.

The optimal PM schedules are presented in Table~\ref{10years}.
All components are planned for PM during months $42$ and~$85$,
    with an optimal total cost for maintenance of $\$4147000$.
Then, when doubling the costs $\PMcost{j}$, the optimal schedule alters to 
    comprise PM of all components only once, in month $66$, 
    with an optimal total cost of $\$5525000$.
%
\begin{table}[!ht]
\newcommand{\tabincell}[2]{\footnotesize \begin{tabular}{@{}#1@{}}#2\end{tabular}}
\centering
{\small
\begin{tabular}{@{}c@{\quad}c@{\quad}c@{\quad}c@{\quad}c@{\quad}c@{\quad}c@{}}
    \hline
{\footnotesize $\PMcost{j}$, $j \in \mathcal{J}$} &
    \multicolumn{4}{c}{\footnotesize Planned PM for components [month]} & 
    {\footnotesize Total cost} & {\footnotesize Computing} \\
{\footnotesize $[\$ 1000]$} & {\footnotesize Rotor} & {\footnotesize Main bearing}
    & {\footnotesize Gearbox} & {\footnotesize Generator} & 
    {\footnotesize $[\$ 1000]$} & {\footnotesize times [s]} \\
    \hline 
    \tabincell{c}{28 \ 15 \ 38 \ 25 \\{}} & 
    \tabincell{c}{42\\85} &
    \tabincell{c}{42\\85} & 
    \tabincell{c}{42\\85} & 
    \tabincell{c}{42\\85} &
    \tabincell{c}{4147\\{}} &
    \tabincell{c}{60+0.53\\{}} \\ 
    \hline
    \tabincell{c}{56 \ 30 \ 76 \ 50} &
    \tabincell{c}{66} &
    \tabincell{c}{66} & 
    \tabincell{c}{66} & 
    \tabincell{c}{66} &
    \tabincell{c}{5525} &
    \tabincell{c}{62+0.24} \\ 
    \hline
\end{tabular}
}
\caption{Optimal schedules (months) for planned PM from the model
    \eqref{model:full-contract}, 
    over a planning period of $10$ years, i.e., $120$ months, for 
    $\mobcostPMf{t} = \$ 50000$ and two levels of
    the costs $\PMcost{j}$, $j \in \mathcal{J}$.
    Computing times represent the calculation of model parameters using Matlab
        ($60$~s and $62$~s, respectively)vand solving the model using CPLEX
        ($0.53$~s and $0.24$~s, respectively)}
    \label{10years}
\end{table}
 
We next compare the model~\eqref{model:full-contract} over ten years with the
    pure CM strategy defined in Section \ref{sec:pureCMcost},
    in which maintenance is performed only when a component fails,
    a {\sl constant interval} (CI) {\sl policy} developed 
    in~\cite{tian2011condition}, and the PMSPIC model developed 
    in~\cite{gustavsson2014preventive}.
As reported in Table~\ref{tab:compare_availibility},
    for the pure CM strategy, the total maintenance cost is $\$6897000$,
    which is $66~\%$ higher than that of our optimal schedule;
for the CI policy, the total cost for maintenance is around $\$5966000$,
    which is $39~\%$ higher than the optimal schedule;
the total maintenance cost resulting from the PMSPIC model is around
    $\$4948000$, which is $20~\%$ higher than that of our optimal schedule. 
\begin{table}[!ht]
\newcommand{\tabincell}[2]{\footnotesize \begin{tabular}{@{}#1@{}}#2\end{tabular}}
\centering
{\small
\begin{tabular}{@{}l@{\quad}c@{\quad}c@{\quad}c@{\quad}c@{}}
    \hline
    \tabincell{c}{Planning method} & 
    \tabincell{c}{Total costs} &
    \tabincell{c}{Relative} &
    \tabincell{c}{Downtime} &
    \tabincell{c}{Availability} \\
    & \tabincell{c}{[$\$ 1000$]} &
    \tabincell{c}{total costs} &
    \tabincell{c}{[months]} &
    \tabincell{c}{[$\%$]} \\
    \hline
    \tabincell{c}{Model \eqref{model:full-contract}} & 
    \tabincell{c}{4147} & 
    \tabincell{c}{1} & 
    \tabincell{c}{1.41} & 
    \tabincell{c}{98.7} \\
    \tabincell{c}{Pure CM strategy} & 
    \tabincell{c}{6897} &
    \tabincell{c}{1.66} & 
    \tabincell{c}{3.75} & 
    \tabincell{c}{96.9} \\
    \tabincell{c}{CI policy} & 
    \tabincell{c}{5762} & 
    \tabincell{c}{1.39} &
    \tabincell{c}{1.54} &
    \tabincell{c}{98.6} \\
    \tabincell{c}{PMSPIC} &
    \tabincell{c}{4948} &
    \tabincell{c}{1.20} &
    \tabincell{c}{1.46} &
    \tabincell{c}{98.7} \\ 
    \hline
\end{tabular}
}
\caption{Total costs ($\$ 1000$) and relative total costs for maintenance, 
        average downtime (month), and technical availability level ($\%$),
        during 10 years for the four maintenance policies,
        with the mobilization costs $\mobcostPMd{t} = \$ 0$ and
        $\mobcostPMf{t} = \$ 50000$.}
\label{tab:compare_availibility}
\end{table}

The resulting technical availability from each of the four methods is also
    presented in Table~\ref{tab:compare_availibility}. 
For the pure CM strategy (with no planned PM), the average number of failures for one wind turbine 
    over 10 years is $3.75$, meaning that the average downtime for each wind turbine is $3.75$ months. 
Applying the model~\eqref{model:full-contract} results in an average number of
    failures per turbine is $1.24$, this is due to the property of the model to result
    in more PM occasions when the probability of a component failure increases.
Combined with two scheduled PM activities for each of the ten turbines, 
    the average downtime for each wind turbine is $1.57$ months, which is less than half of the downtime
    under the pure CM strategy. 

The technical availability resulting from the model~\eqref{model:full-contract}
    is thus considerably higher than that of the pure CM strategy,
    while also generating a substantially lower cost.
Comparing with the CI policy and the PMSPIC model, the technical availability
    is approximately the same, but the total maintenance cost resulting
    from either of these two methods is considerably higher.
From the comparison of these four methods, we conclude that the 
    model~\eqref{model:full-contract} yields the highest technical availability
    at the lowest total cost for maintenance during ten years.

Nowadays, it is common to have a $97~\%$ production-based availability guarantee. 
With the three PM methods, the expected level of availability is around $98.7~\%$,
    while with the pure CM strategy, the expected level of availability is below $97~\%$. 
Hence, we conclude that planning for PM will provide great benefits
 by its ability to provide a substantially increased technical availability level.


\subsection{Case study 2: Lifetime planning with no insurance contract and
    sudden component failures} 
 
We next consider lifetime planning with sudden component failures.
The failures are sampled from the Weibull distribution, as defined in
    Table~\ref{tab:Task1_parameters}.
    
When a component fails, it has to undergo CM while the PM plan needs
    to be rescheduled from the point in time of the failure.
The new schedule is then used until another component fails, at which 
    time the PM plan is rescheduled. 
This process is repeated until the end of life of the wind farm.

We first plan PM for the ten turbines over a ten years maintenance period
    in phase~\ref{normal-phase}. 
Hence, the model~\eqref{model:normal} is applied; the resulting PM schedule is
    presented in Table~\ref{tab:normal-plan}.
%
\begin{table}[!ht]
\newcommand{\tabincell}[2]{\footnotesize \begin{tabular}{@{}#1@{}}#2\end{tabular}}
\centering
{\small
\begin{tabular}{@{}c@{\quad}c@{\quad}c@{\quad}c@{\quad}c@{\quad}c@{\quad}c@{}}
    \hline
    {\footnotesize Turbine} &
        \multicolumn{4}{c}{\footnotesize Planned PM for components [month]} &
    {\footnotesize Cost} & {\footnotesize Computing} \\
    \# & {\footnotesize Rotor} & {\footnotesize Main bearing} & {\footnotesize Gearbox} & {\footnotesize Generator} & {\footnotesize $[\$ 1000]$} & {\footnotesize times [s]} \\
    \hline 
    \tabincell{c}{1--10\\{}\\{}} &
    \tabincell{c}{39\\78\\120} & 
    \tabincell{c}{39\\78\\120} & 
    \tabincell{c}{39\\78\\120} & 
    \tabincell{c}{39\\78\\120} &
    \tabincell{c}{5294/44.1\\{}\\{}} &
    \tabincell{c}{28+0.6\\{}\\{}} \\ 
    \hline
\end{tabular}
}
\caption{Optimal schedule (months) for planned PM from the model
    \eqref{model:normal}, 
    over a planning period of $10$ years, i.e., during the months $1$--$120$.
    'Cost' denotes the "total expected cost for performed and planned maintenance
    from month $1$"/"average expected cost per month for performed and planned
    maintenance".
    'Computing times' represent calculation of model parameters using Matlab ($28$~s)
    and solving the model using CPLEX ($0.6$~s)}
  \label{tab:normal-plan}
\end{table}

A gearbox break-down in turbine~$1$ is then sampled at month~$31$.
The new PM schedule, calculated using the model \eqref{model:normal},
    is presented in 
    Table~\ref{tab:normal-plan-1}.\footnote{The {\bf bold} figures in Tables~\ref{tab:normal-plan-1}--\ref{tab:normal-plan-4} denote CM occasions.}
The total cost denotes the sum of costs for all PM and CM performed before
    the planning period, plus the cost of PM during the planning period at hand.
Since the planned PM of gearbox~$1$ in month~$39$ is replaced by CM in month~$31$,
    and after that all components receive PM three times during the planning
    period,
    the total cost of the new schedule is larger than that of the previously
    planned schedule (before the gearbox break-down).
The total cost per month is, however, lower for the new schedule
    (i.e., $\$6029000/150 \approx \$40200$) than for the previous schedule
    (i.e., $\$5294000/120 \approx \$44100$).
%
\begin{table}[!ht]
\newcommand{\tabincell}[2]{\footnotesize \begin{tabular}{@{}#1@{}}#2\end{tabular}}
\centering
{\small
\begin{tabular}{@{}c@{\quad}c@{\quad}c@{\quad}c@{\quad}c@{\quad}c@{\quad}c@{}}
    \hline
    {\footnotesize Turbine} & \multicolumn{4}{c}{\footnotesize Planned PM for components [month]} &
    {\footnotesize Cost} & {\footnotesize Computing} \\
    \# & {\footnotesize Rotor} & {\footnotesize Main bearing} & {\footnotesize Gearbox} & {\footnotesize Generator} & {\footnotesize $[\$ 1000]$} & {\footnotesize times [s]} \\
    \hline 
    \tabincell{c}{1\\{}\\{}\\{}\\{}\\{}} & 
    \tabincell{c}{{}\\{\sl 54}\\{}\\102\\{}\\150} & 
    \tabincell{c}{{}\\{\sl 54}\\{}\\102\\{}\\150} & 
    \tabincell{c}{{\bf 31}\\{}\\73\\{}\\111\\150} & 
    \tabincell{c}{{}\\{\sl 54}\\{}\\102\\{}\\150} &
    \tabincell{c}{6029/40.2\\{}\\{}\\{}\\{}\\{}} &
    \tabincell{c}{98+18.6\\{}\\{}\\{}\\{}\\{}}\\ 
    \cline{1-5}
    \tabincell{l}{2--10\\{}\\{}\\{}\\{}\\{}} &
    \tabincell{c}{{}\\{\sl 54}\\{}\\102\\{}\\150} & 
    \tabincell{c}{{}\\{\sl 54}\\{}\\102\\{}\\150} & 
    \tabincell{c}{{\sl 39}\\{}\\73\\{}\\111\\150} & 
    \tabincell{c}{{}\\{\sl 54}\\{}\\102\\{}\\150} &
    \tabincell{c}{\\{}\\{}} \\ 
    \hline
\end{tabular}
}
\caption{Optimal schedules (months) for planned PM from the model
    \eqref{model:normal}---with a sudden break of the gearbox
    of turbine~$1$ in month~$31$---over a planning period of~$10$ years,
    i.e., during the months~$31$--$150$.
    'Cost' denotes the "total expected cost for performed and planned maintenance
    from month $1$"/"average expected cost per month for performed and planned
    maintenance".
    'Computing times': calculating model parameters using Matlab ($98$~s);
    solving the model using CPLEX ($18.6$~s)}
  \label{tab:normal-plan-1}
\end{table}

Then, a break-down is sampled of a main bearing in turbine~$2$ in month~$71$.
The PM plan in Table~\ref{tab:normal-plan-1} is thus followed until all the
    planned replacements in month $54$ are
    performed.\footnote{The {\sl slanted} figures in
    Tables~\ref{tab:normal-plan-1}--\ref{tab:normal-plan-4} denote the planned PM
    occasions that are {\sl actually performed}; note that this information
    is not revealed until the "next" component break-down occurs.}
Then, a new schedule is calculated from month~$71$;
    it is presented in Table~\ref{tab:normal-plan-2}.
The total cost $\$ 6655000$ equals the sum of costs for CM and PM already performed
    (i.e., CM of the gearbox in turbine $1$ in month $31$, and the performed PM
    in months $39$ and $54$
    see Table~\ref{tab:normal-plan-1}), CM of the main bearing
    in month $71$, and the expected cost for all the planned PM
    (see Table~\ref{tab:normal-plan-2}).
%
\begin{table}[!ht]
\newcommand{\tabincell}[2]{\footnotesize \begin{tabular}{@{}#1@{}}#2\end{tabular}}
\centering
{\small
\begin{tabular}{@{}c@{\quad}c@{\quad}c@{\quad}c@{\quad}c@{\quad}c@{\quad}c@{}}
    \hline
    {\footnotesize Turbine} & \multicolumn{4}{c}{\footnotesize Planned PM for components [month]} &
    {\footnotesize Cost} & {\footnotesize Computing} \\
    \# & {\footnotesize Rotor} & {\footnotesize Main bearing} & {\footnotesize Gearbox} & {\footnotesize Generator} & {\footnotesize $[\$ 1000]$} & {\footnotesize times [s]} \\
    \hline 
    \tabincell{c}{2\\{}\\{}\\{}\\{}} & 
    \tabincell{c}{{}\\{\sl 99}\\{}\\145\\189} & 
    \tabincell{c}{{\bf 71}\\{}\\{\sl 110}\\145\\189} & 
    \tabincell{c}{{\sl 71}\\{}\\{\sl 110}\\145\\189} & 
    \tabincell{c}{{}\\{\sl 99}\\{}\\145\\189} &
    \tabincell{c}{6655/35.0\\{}\\{}\\{}\\{}} &
    \tabincell{c}{171+72\\{}\\{}\\{}\\{}} \\ 
    \cline{1-5}
    \tabincell{l}{1, 3--10\\{}\\{}\\{}\\{}} &
    \tabincell{c}{{}\\{\sl 99}\\{}\\145\\189} & 
    \tabincell{c}{{}\\{\sl 99}\\{}\\145\\189} & 
    \tabincell{c}{{\sl 71}\\{}\\{\sl 110}\\145\\189} & 
    \tabincell{c}{{}\\{\sl 99}\\{}\\145\\189} &
        \tabincell{c}{\\{}\\{}} \\ 
    \hline
\end{tabular}
}
\caption{Optimal schedules (months) for planned PM from the model
    \eqref{model:normal}---with a sudden break of the main bearing of
    turbine~$2$ in month~$71$---over~$10$ years, i.e., months~$71$--$190$.
    'Cost' denotes the "total expected cost for performed and planned maintenance
    from month $1$"/"average expected cost per month for performed and planned
    maintenance".
    'Computing times': calculating model parameters using Matlab ($171$~s);
    solving the model using CPLEX ($72$~s)}
\label{tab:normal-plan-2}
\end{table}

Then, a break-down is sampled of a rotor in turbine~$3$ in month~$131$.
The PM plan in Table~\ref{tab:normal-plan-2} is thus followed until all the
    planned replacements in month~$110$ are performed. 
Then, a new schedule is calculated from month~$131$ until the end of life
    of the wind farm, hence using the model~\eqref{eq:hatFmodel}.
The resulting PM schedule is presented in Table~\ref{tab:normal-plan-3}.
The total cost $\$7107000$ is the sum of costs for CM and PM already performed
    (i.e., CM of the gearbox in turbine~$1$ in month~$31$,
    CM of the main bearing of turbine~$2$ in month~$71$,
    and the performed PM in months $39$, $54$, $71$, $99$, and~$110$,
    according to Tables~\ref{tab:normal-plan-1}--\ref{tab:normal-plan-2}),
    CM of the rotor in month~$131$, and the expected cost for all the planned PM
    in months $150$ and~$195$ (according to Table~\ref{tab:normal-plan-3}).
%
\begin{table}[!ht]
\newcommand{\tabincell}[2]{\footnotesize \begin{tabular}{@{}#1@{}}#2\end{tabular}}
\centering
{\small
\begin{tabular}{@{}c@{\quad}c@{\quad}c@{\quad}c@{\quad}c@{\quad}c@{\quad}c@{}}
    \hline
    {\footnotesize Turbine} & \multicolumn{4}{c}{\footnotesize Planned PM for components [month]} &
    {\footnotesize Cost} & {\footnotesize Computing} \\
    \# & {\footnotesize Rotor} & {\footnotesize Main bearing} & {\footnotesize Gearbox} & {\footnotesize Generator} & {\footnotesize $[\$ 1000]$} & {\footnotesize times [s]} \\
    \hline 
    \tabincell{c}{3\\{}\\{}} & 
    \tabincell{c}{{\bf 131}\\{}\\195} & 
    \tabincell{c}{{}\\{\sl 150}\\195} & 
    \tabincell{c}{{}\\{\sl 150}\\195} & 
    \tabincell{c}{{}\\{\sl 150}\\195} &
    \tabincell{c}{7107/29.6\\{}\\{}} &
    \tabincell{c}{176+98\\{}\\{}} \\ 
    \cline{1-5}
    \tabincell{l}{1--2, 4--10\\{}} &
    \tabincell{c}{{\sl 150}\\195} & 
    \tabincell{c}{{\sl 150}\\195} & 
    \tabincell{c}{{\sl 150}\\195} & 
    \tabincell{c}{{\sl 150}\\195} &
    \tabincell{c}{{}\\{}} \\ 
    \hline
\end{tabular}
}
\caption{Optimal schedules (months) for planned PM from the model
    \eqref{eq:hatFmodel}---with a sudden break of a rotor of turbine~$3$
    in month~$131$---for a planning period comprising the rest of the
    wind farm's life, i.e., during months~$131$--$240$.
    'Cost' denotes the "total expected cost for performed and planned maintenance
    from month $1$"/"average expected cost per month for performed and planned
    maintenance".
    'Computing times': calculating model parameters using Matlab ($176$~s);
    solving the model using CPLEX ($98$~s)}
\label{tab:normal-plan-3}
\end{table}

Then, a break-down is sampled of the generator of turbine~$4$ in month~$181$.
The PM plan in Table~\ref{tab:normal-plan-3} is thus followed until all the
    planned replacements in month~$150$ are performed. 
Then, a new schedule is calculated from month~$181$;
    it is presented in Table~\ref{tab:normal-plan-4}.
The total cost $\$6463000$ is the sum of costs for CM and PM already performed
    (i.e., CM of the gearbox in turbine~$1$ in month~$31$,
    CM of the main bearing of turbine~$2$ in month~$71$,
    CM of the rotor of turbine~$3$ in month~$131$,
    and the performed PM in months $39$, $54$, $71$, $99$, $110$, and~$150$,
    according to Tables~\ref{tab:normal-plan-1}--\ref{tab:normal-plan-3}),
    CM of the generator in month~$181$,
    and all the planned PM in month~$195$ (according to Table~\ref{tab:normal-plan-4}).
%
\begin{table}[!ht]
\newcommand{\tabincell}[2]{\footnotesize \begin{tabular}{@{}#1@{}}#2\end{tabular}}
\centering
{\small
\begin{tabular}{@{}c@{\quad}c@{\quad}c@{\quad}c@{\quad}c@{\quad}c@{\quad}c@{}}
    \hline
    {\footnotesize Turbine} &
        \multicolumn{4}{c@{}}{\footnotesize Planned PM for components [month]} &
        {\footnotesize Cost} & {\footnotesize Computing} \\
    \# & {\footnotesize Rotor} & {\footnotesize Main bearing} &
        {\footnotesize Gearbox} & {\footnotesize Generator} &
        {\footnotesize $[\$ 1000]$} & {\footnotesize times [s]} \\
    \hline 
    \tabincell{c}{4\\{}} & 
    \tabincell{c}{{}\\{\sl 195}} & 
    \tabincell{c}{{}\\{\sl 195}} & 
    \tabincell{c}{{}\\{\sl 195}} & 
    \tabincell{c}{{\bf 181}\\{}} &
    \tabincell{c}{6463/26.9\\{}} &
    \tabincell{c}{86+15\\{}} \\ 
    \cline{1-5}
    \tabincell{l}{1--3, 5--10} &
    \tabincell{c}{{\sl 195}} & 
    \tabincell{c}{{\sl 195}} & 
    \tabincell{c}{{\sl 195}} & 
    \tabincell{c}{{\sl 195}} &
    \tabincell{c}{} \\ 
    \hline
\end{tabular}
}
\caption{Optimal schedules (months) for planned PM from the model
    \eqref{eq:hatFmodel}---with a sudden break of the generator of turbine~$4$
    in month~$181$---for a planning period comprising the rest of the wind 
    farm's life, i.e., during months~$181$--$240$.
    'Cost' denotes the "total expected cost for performed and planned maintenance
    from month $1$"/"average expected cost per month for performed and planned
    maintenance".
    'Computing times': calculating model parameters using Matlab ($86$~s);
    solving the model using CPLEX ($15$~s)}
  \label{tab:normal-plan-4}
\end{table}

Then, a break-down of the gearbox of turbine~$5$ is sampled in month~$235$.
The PM plan in Table~\ref{tab:normal-plan-4} is thus followed until all the
    planned replacements in month~$195$ are performed. 
Since it is close to the end of life of the wind turbine, by comparing the
    production loss for one turbine during the five remaining months and
    the cost of CM, we conclude that it will not be beneficial to
    maintain this gearbox during month~$235$.
 
The final total maintenance cost, including the downtime cost resulting from 
    the updated PM plans at CM occasions, as well as the downtime cost of
    one turbine after the breakdown in month~$235$, equals $\$5798000$,
    which corresponds to an average of $\$24200$ per month during the life
    of the wind farm. 
The expected number of breakdowns during this wind farm's life---when
    following the PM plan---is larger than~$20$.
Due to not experiencing as many breakdowns as expected,
    the final total maintenance cost is thus lower than expected---compare
    with the numbers in Tables~\ref{tab:normal-plan-3}--\ref{tab:normal-plan-4}.

\subsection{Case study 3: Lifetime planning with no insurance contract} \label{sec:lifetimeplanning} 

We next consider a case with no insurance contract,
    i.e., a contract of type~\ref{no-contract}, over~$20$ years,
    ($240$ months) during the phase~\ref{end-of-life-phase}.
Hence, the model~\eqref{eq:hatFmodel} is applied;
    the resulting PM schedule is presented in Table~\ref{tab:whole-life}.
%
 \begin{table}[!ht]
\newcommand{\tabincell}[2]{\footnotesize \begin{tabular}{@{}#1@{}}#2\end{tabular}}
\centering
{\small
\begin{tabular}{@{}c@{\quad}c@{\quad}c@{\quad}c@{\quad}c@{\quad}c@{}}
    \hline
    \multicolumn{4}{c}{\footnotesize Planned PM for components [month]} &
        {\footnotesize Total cost} & {\footnotesize Computing} \\
    {\footnotesize Rotor} & {\footnotesize Main bearing} &
        {\footnotesize Gearbox} & {\footnotesize Generator} &
        {\footnotesize $[\$ 1000]$} & {\footnotesize times [s]} \\
    \hline 
    \tabincell{c}{{}\\49\\{}\\99\\{}\\148\\{}\\198} & 
    \tabincell{c}{{}\\49\\{}\\99\\{}\\148\\{}\\198} & 
    \tabincell{c}{39\\{}\\78\\{}\\119\\{}\\159\\198} & 
    \tabincell{c}{{}\\49\\{}\\99\\{}\\148\\{}\\198} &
    \tabincell{c}{9153\\{}\\{}\\{}\\{}\\{}\\{}\\{}} &
    \tabincell{c}{213+3.4\\{}\\{}\\{}\\{}\\{}\\{}\\{}} \\ 
    \hline
\end{tabular}
}
\caption{Optimal schedules (months) for planned PM from the model
    \eqref{eq:hatFmodel}, 
    over a planning period of~$20$ years, i.e.,~$240$ months.
    Computing times represent calculation of model parameters using Matlab
    ($213$~s) and solving the model using CPLEX ($3.4$~s)}
  \label{tab:whole-life}
\end{table}
The total maintenance cost over these~$20$ years is $\$9153000$,
    which is larger than twice the maintenance cost,
    i.e., $\$4147000$ over ten years, as computed in Case study~$1$
    considering the phase~\ref{end-of-full-contract-phase}, the end of
    a full service production-based maintenance contract.
Hence, at the end of phase~\ref{end-of-full-contract-phase}, even though the
    components are old, there is no additional cost for compensating age.
Moreover, the total cost over $20$ years in phase~\ref{end-of-life-phase} is 
    lower than twice the cost $\$5294000$ of the first ten years
    plan in Case study~$2$, considering the phase~\ref{normal-phase}.
This is due to the fact that at the end of life, the component ages need
    not be considered. 

We then look further into the computational complexity of our models by
    redefining the time steps to comprise three days instead of one month,
    resulting in the number of time steps in the model being multiplied by ten;
    the resulting PM schedule is presented in Table~\ref{tab:3days}.
%
\begin{table}[!htb]
\newcommand{\tabincell}[2]{\footnotesize \begin{tabular}{@{}#1@{}}#2\end{tabular}}
\centering
{\small
\begin{tabular}{@{}c@{\quad}c@{\quad}c@{\quad}c@{\quad}c@{\quad}c@{}}
    \hline
 \multicolumn{4}{c}{\footnotesize Planned PM for components [3 days]} &
        {\footnotesize Total cost} & {\footnotesize Computing} \\
    {\footnotesize Rotor} & {\footnotesize Main bearing} &
        {\footnotesize Gearbox} & {\footnotesize Generator} &
        {\footnotesize $[\$ 1000]$} & {\footnotesize times [min]} \\
    \hline 
    \tabincell{c}{{}\\495\\{}\\990\\{}\\1486\\{}\\1982} & 
    \tabincell{c}{{}\\495\\{}\\990\\{}\\1486\\{}\\1982} & 
    \tabincell{c}{394\\{}\\788\\{}\\1190\\{}\\1586\\1982} & 
    \tabincell{c}{{}\\495\\{}\\990\\{}\\1486\\{}\\1982} &
    \tabincell{c}{9129\\{}\\{}\\{}\\{}\\{}\\{}\\{}} &
    \tabincell{c}{65+319\\{}\\{}\\{}\\{}\\{}\\{}\\{}} \\ 
    \hline
\end{tabular}
}
\caption{Optimal schedules ($3$~days) for planned PM from the model
    \eqref{eq:hatFmodel}, 
    over a planning period of~$20$ years, i.e.,~$240$ months.
    Computing times represent calculation of model parameters using Matlab
        ($65$ min) and solving the model using CPLEX ($319$ min)}
  \label{tab:3days}
\end{table}
Comparing with the results from employing time steps of one month,
    the schedules are refined, but essentially equivalent.
Hence, the finer time discretization mainly lead to substantially increased
    computing times while the results are not improved.
We conclude that the time steps of one month should be kept.

\subsection{Case study~$4$: Lifetime planning with a basic insurance contract}

We next investigate the case when---at a component failure---with the
    probability $\delta \in [0,1]$ the insurance company will pay for
    a new component. 
Then the cost of a new component plus the logistics cost for the replacement,
    i.e., $\CMCOST{j}$, is reduced by the probability times the component cost,
    i.e., $\delta \CMcompCOST{j}$, according to the expressions in~\eqref{eq:Cdot}.
We employ the values of $\CMcompCOST{j}$ from~\citet[Table~2]{tian2011condition}
    and which are given by $\$112000$, $\$60000$, $\$152000$, and $\$100000$
    for the rotor, the main bearing, the gearbox, and the generator, respectively.

Here we apply the 
modified version of model~\eqref{eq:hatFmodel}, with the cost defined 
    in~\eqref{eq:interval-cost-vC} altered analogously as the costs in the
    formula~\eqref{eq:Cdot}.
For different values of
    $\delta \in \{ 1.0, 0.8, 0.6, 0.4, 0.2 \}$, the optimization model yields the corresponding
    optimal PM schedules presented in Table~\ref{insurance}.
Note that for $\delta \in \{ 0.8, 0.6 \}$ the optimal schedules are equal
    and that also $\delta \in \{ 0.4, 0.2 \}$ yield equal optimal schedules.
Note also that for large values of $\delta$---representing a high probability
    that the insurance company will pay for the component cost at
    a failure---the wind farm owner should plan for fewer PM activities,
    i.e., the number of planned PM occasions is non-increasing with an
    increasing value of $\delta$.
\begin{table}[!htb]
\newcommand{\tabincell}[2]{\footnotesize \begin{tabular}{@{}#1@{}}#2\end{tabular}}
\centering
{\small
\begin{tabular}{@{}c@{\quad}c@{\quad}c@{\quad}c@{\quad}c@{\quad}c@{}}
    \hline
    {\footnotesize $\delta$} &
    \multicolumn{4}{c}{\footnotesize Planned PM for components [month]} & 
    {\footnotesize $F^*(\delta)$} \\
    & {\footnotesize Rotor} & {\footnotesize Main bearing} &
    {\footnotesize Gearbox} & {\footnotesize Generator} &
    {\footnotesize $[\$ 1000]$} \\
    \hline 
 \tabincell{c}{1.0\\{}} &
    \tabincell{c}{{90}\\181} &  
    \tabincell{c}{{90}\\181} & 
    \tabincell{c}{90\\{181}} & 
    \tabincell{c}{{90}\\{181}} &
    \tabincell{c}{5292\\{}} \\ 
    \hline
\tabincell{l}{0.8\\{}\\{}} & 
    \tabincell{c}{63\\123\\186} &  
    \tabincell{c}{63\\123\\186} &  
    \tabincell{c}{63\\123\\186} &  
    \tabincell{c}{63\\123\\186} &
    \tabincell{c}{6471\\{}\\{}} \\ 
    \hline
    \tabincell{l}{0.6\\{}\\{}} & 
    \tabincell{c}{63\\123\\186} &  
    \tabincell{c}{63\\123\\186} &  
    \tabincell{c}{63\\123\\186} &  
    \tabincell{c}{63\\123\\186} &
    \tabincell{c}{7385\\{}\\{}} \\ 
    \hline
    \tabincell{l}{0.4\\{}\\{}\\{}} & 
    \tabincell{c}{49\\99\\148\\198} &  
    \tabincell{c}{49\\99\\148\\198} &  
    \tabincell{c}{49\\99\\148\\198} &  
    \tabincell{c}{49\\99\\148\\198} &
    \tabincell{c}{8054\\{}\\{}\\{}} \\ 
    \hline
    \tabincell{l}{0.2\\{}\\{}\\{}} & 
    \tabincell{c}{49\\99\\148\\198} &  
    \tabincell{c}{49\\99\\148\\198} &  
    \tabincell{c}{49\\99\\148\\198} &  
    \tabincell{c}{49\\99\\148\\198} &
    \tabincell{c}{8620\\{}\\{}\\{}} \\ 
    \hline
\end{tabular}
}
\caption{Optimal schedules (months) for planned PM resulting from the 
    modified version of model~\eqref{eq:hatFmodel},
    for different values of $\delta \in [0,1]$,
    over a planning period of~$20$ years, i.e.,~$240$ months. 
    $F^*(\delta)$ denotes the total cost of maintenance during the~$20$ years
    under a $\delta$-insurance contract}
\label{insurance}
\end{table}
The optimal value of the lifetime PM schedule for the case when there is
    no insurance contract (i.e., for $\delta = 0.0$), as presented in 
    Table~\ref{tab:whole-life}, equals $F^*(0) = \$9153000$.
    
We conclude that the benefit of planned PM depends on the terms of the insurance,
    and that a reasonable price of a $\delta$-insurance maintenance contract can be computed as the difference $F^*(0)\!-\!F^*(\delta)$.

\section{Conclusion} \label{sec:conclusion}

\noindent
This paper models the optimal planning of preventive maintenance (PM)
    for a wind farm.
Our models take into account component failures and necessary corrective
    maintenance (CM) as well as the effect of the expected costs for CM
    on the PM planning.
We provide optimization models for different kinds of contracts and during
    different phases of the life of the wind farm.

Case study~$1$ reveals that, as applied to a full service production-based
    maintenance contract, a comparison of our model with a pure CM strategy,
    a constant interval policy, and a pure PM model, indicates that
    planning for PM while taking into account costs for expected CM
    provides a high technical availability level at a comparatively
    low cost for maintenance.

Case study~$2$ concerns lifetime planning with no insurance contract,
    with a rescheduling of the PM plan each time a component breaks and 
    CM is necessary, throughout the life of the wind farm. 
With only five sampled breakdowns, the eventual total cost for maintenance
    is lower than expected. 
Due to the low PM cost, when one component in one wind turbine breaks,
    it affects the PM plan only for that particular wind turbine, while the
    other turbines' PM plans are mostly unaffected.


Case study~$3$ reveals that when scheduling the maintenance during the 
    entire life of the wind farm, the first PM occasions are almost equidistant
    while the intervals between the later PM occasions are shorter.
This is explained in Section~\ref{parameter_analysis}: 
    for an increasing interval length
    the interval cost between (d)~a PM occasion and the end of life,
    increases faster than that between (a)~two PM occasions.

Since our optimization model is NP-hard (see~\cite{gustavsson2014preventive}),
    the computing time is expected to increase exponentially with the problem
    size  (in terms of numbers of variables and constraints).
This property is confirmed in Case study~$3$, where we investigate a finer
    (by a factor of ten) discretization of time.



Case study~$4$ regards an insurance contract, under which an insurance company
    pays for part of the CM costs caused by (sudden) component failures; 
    the insurance level is then modelled by a probability that the CM
    costs are paid by the insurance company.
Our results show that the number of planned PM occasions for the wind farm
    decreases with an increasing insurance level.
Further, the optimal costs for PM can be used to advise the wind farm owners
    what is a fair price of a certain level of an insurance contract.
%

The optimization models presented in this article regards PM scheduling and
    CM strategies for the stakeholders involved, i.e., wind farm owners,
    maintenance companies, and insurance companies.
Solutions from these models can be used to
    guide on reasonable costs and prices for different types of contracts between the
    stakeholders during different phases of the life of a wind farm.

\section*{Acknowledgements}
\noindent We acknowledge the financial support from the Swedish Wind Power Technology Centre at Chalmers University of Technology.

\section*{References}
\noindent 
\bibliographystyle{plainnat}
\bibliography{main}
\end{document}